\documentclass[11pt]{article}
\usepackage{imports}
\usepackage[utf8]{inputenc}
\usepackage{multicol}
\usepackage{amsmath}
\usepackage{amsfonts}
\usepackage{amssymb}
\usepackage{csquotes}
\usepackage{mathtools}
\usepackage[letterpaper, portrait, margin=1in,bottom=0.5in]{geometry}
\usepackage{graphicx}
    \graphicspath{ {./images/} } 

\usepackage{pst-node}
\usepackage{auto-pst-pdf}
\usepackage{tikz-cd}
\usepackage{enumitem}
\newcommand{\defeq}{\vcentcolon=}

\newcommand{\site}[2]{\cite[#1]{#2}}
\setlist[enumerate]{noitemsep}
\usepackage[normalem]{ulem}

\numberwithin{equation}{section}

\usepackage[numbers]{natbib}

\newenvironment{Proof}{\noindent \textit{Proof.}}{\hfill$\square$}

\newcounter{count}[section]
\counterwithin{count}{section}

\newenvironment{Theorem}[1][]{\refstepcounter{count}\noindent\textbf{Theorem \thecount}#1\textbf{.}}

\newenvironment{definition}[1][\textbf{. }]{\refstepcounter{count}\noindent\textbf{Definition \thecount}#1}

\newenvironment{Lemma}[1][\textbf{. }]{\refstepcounter{count}\noindent\textbf{Lemma \thecount}#1}

\newenvironment{example}[1][\textbf{. }]{\refstepcounter{count}\noindent\textbf{Example \thecount}#1}

\newenvironment{Corollary}[1][\textbf{. }]{\refstepcounter{count}\noindent\textbf{Corollary \thecount}#1}

\newenvironment{remark}[1][\textbf{. }]{\refstepcounter{count}\noindent\textbf{Remark \thecount}#1}

\definecolor{DarkGreen}{RGB}{34,139,34}

\title{New Applications of Ergodic Theory to Sets of Differences}
\author{Kabir Belgikar, Vitaly Bergelson, Gabriel Black, David Kruzel}
\date{}

\begin{document}
\maketitle
\noindent \textbf{Abstract.} We apply the methods of ergodic theory to both simplify and significantly extend some classical results due to Stewart, Tijdeman, and Ruzsa. One of the notable features of our approach is the utilization of pointwise ergodic theory.

\section{Introduction}

Given a set $S \subseteq \N=\{1,2,3,\ldots\}$, the \textit{upper density} of $S$, denoted $\overline{d}(S)$, is defined by
        \begin{equation*}
            \overline{d}(S) = \limsup_{N \rightarrow \infty} \frac{|S \cap \{1,\ldots,N\}|}{N},
        \end{equation*}
and the \textit{lower density} of $S$, denoted $\underline{d}(S)$, is defined by
        \begin{equation*}
            \underline{d}(S) = \liminf_{N \rightarrow \infty} \frac{|S \cap \{1,\ldots,N\}|}{N}.
        \end{equation*}    
If $\overline{d}(S) = \underline{d}(S)$, then the common value is denoted by $d(S)$ and is called the \textit{natural density} of $S$. Following Rusza, we introduce three types of difference sets; note that these sets are symmetric subsets of $\Z$. 
        \begin{equation*}
            \Delta_1(S) = \{n \in \Z: S \cap (S-n) \neq \emptyset\},\footnote{For $S\subseteq \N$ we define $S-n=\{m\in \N: m+n\in S\}$.} \quad \quad \Delta_2(S) = \{n \in \Z: \overline{d}(S \cap (S-n)) > 0 \}, 
        \end{equation*}
        \begin{equation*}
            \Delta_3(S) = \{n \in \Z: |S \cap (S-n) | = \infty \}. \footnote{Our definitions of $\Delta_2(S)$ and $\Delta_3(S)$ are switched from how Ruzsa defines them in \cite{Ruz78}, on account of the fact that we will not deal with sets of the form $\Delta_3(S)$ outside of this introduction.  }
        \end{equation*}
    \indent Suppose that $S_1,\ldots,S_k \sseq \N$ and $\overline{d}(S_i) > 0$ for $1 \leq i \leq k$. In \cite{ST79}, Stewart and Tijdeman answered in the affirmative a question of Erdős by showing 
    that the set $C\defeq\bigcap_{i=1}^k \Delta_3(S_i)$ has positive lower density. Moreover, they showed that $C$ is \textit{syndetic}, meaning that there are integers $m_1,\ldots,m_\ell$ such that
        \begin{equation}
            \bigcup_{i=1}^\ell \left( C + m_i \right) = \Z.\footnote{A set $C\subseteq \N$ is said to be syndetic if there exist $m_1,\dots,m_{\ell}\in \N$ such that $\bigcup_{i=1}^\ell (C-m_i)=\N$.}
        \end{equation}
    Stewart and Tijdeman's result was amplified by Rusza as follows.
    
    \begin{Theorem}[ (\site{Theorem 1 and Theorem 2}{Ruz78})] \label{Ruz}
        If $S_1,\ldots, S_k \subseteq \N$ satisfy $\overline{d}(S_i) > 0$ for each $i=1,\dots,k$, then there exists $S \subseteq \N$ such that $d(S) \geq \prod_{i=1}^k \overline{d}(S_i)$ and $\Delta_1 (S) \subseteq D \defeq \bigcap_{i=1}^k \Delta_2 (S_i)$. Furthermore, there are $m_1,\ldots,m_\ell \in \Z$, where
            \begin{equation}
            \ell \leq \prod_{i=1}^k 1 / \overline{d}(S_i), 
            \end{equation}
        such that $\bigcup_{i=1}^\ell (D + m_i) = \Z$. 
    \end{Theorem}

     The goal of this paper is to offer an ergodic approach to Theorem \ref{Ruz} which not only provides an alternative proof, but also leads to various generalizations which range from establishing versions of Theorem \ref{Ruz} in the setup of countably infinite \textit{amenable cancellative} semigroups\footnote{Amenable cancellative semigroups form a natural framework where notions of density can be defined. A precise definition of amenability is given in subsection \ref{Sec 4.2}. 
     } to variants of Theorem \ref{Ruz} which involve modified versions of the sets $\Delta_1(S)$ and $\Delta_2(S)$ such as
    \begin{equation}\label{1.5}\Delta_{1,c}(S)=\{n\in \N:S\cap (S-[n^c])\neq \emptyset\}\quad \text{ and }\quad \Delta_{2,c}(S)=\{n\in \N:\ud(S\cap (S-[ n^c]))>0\}\end{equation}
    with $c > 0$, or, more generally,
     \begin{equation}\label{1.6}\Delta_{1,g}(S)=\{n\in \N:S\cap (S-[g(n)])\neq \emptyset\}\quad \text{ and }\quad \Delta_{2,g}(S)=\{n\in \N:\ud(S\cap (S-[ g(n)]))>0\}\end{equation}
    for certain functions $g:\N \rightarrow [1,\infty)$ (in formulas (\ref{1.5}) and (\ref{1.6}), $[\cdot]$ denotes the integral part). An interesting new variant of Theorem \ref{Ruz} is provided by the following theorem, which deals with the function $g(n) = n^c$ and is a special case of the much more general Theorem \ref{3.9}, to be proved in subsection \ref{Sec 3.2}.


\begin{Theorem} \label{1.2}
        Let $c>0$, $c\not\in \N$. If $S_1,\dots,S_k\subseteq \N$ satisfy $\ud(S_i)>0$ for each $i=1,\dots,k$, then there is $S\subseteq \N$ with $d(S)\geq \prod_{i=1}^k \ud(S_i)$ and $ \Delta_{1,c} (S)\subseteq D_c\defeq \bigcap_{i=1}^k \Delta_{2,c} (S_i)$.
        Furthermore, we have that $\ud(\Delta_{1,c} (S))\geq \prod_{i=1}^k \ud(S_i)$. In addition, the set $D_c$ is thick, meaning that for all $n\in\N$ there is an $a\in\N$ such that $\{a,a+1,\dots,a+n\}\subseteq D_c$.
    \end{Theorem}

    The reader may wonder why we require $c \notin \N$ in the statement of Theorem \ref{1.2}. The reason for this (to be revealed in subsection \ref{ergodic sequences}) is that for $c \notin \N$, the sequence $([n^c])_{n=1}^\infty$ is \textit{pointwise ergodic} (see Definition \ref{ergodic defs}), which in turn leads to the rather sharp inequality $d(S) \geq \prod_{i=1}^k \overline{d}(S_i)$ appearing in the statement of Theorem \ref{1.2}. However, even when $c \in \N$, we still have the following nontrivial statement which naturally complements Theorem \ref{1.2} and is a special case of Theorem \ref{3.13}.
    
    \begin{Theorem} \label{c in N}
        Let $c \in \N$. If $S_1,\ldots,S_k \sseq \N$ satisfy $\overline{d}(S_i) > 0$ for each $i = 1,\ldots,k$, then there is $S \sseq \N$ with $d(S) > 0$, and $\Delta_{1,c}(S) \sseq  D_c \defeq \bigcap_{i=1}^k \Delta_{2,c} (S_i)$. Furthermore, the set $D_c$ is syndetic. 
    \end{Theorem}

    Theorem \ref{c in N} is weaker than Theorem \ref{1.2} in the sense that it does not provide a bound for $d(S)$. But, unlike Theorem \ref{1.2}, Theorem \ref{c in N} guarantees syndeticity of the set $D_c$. 
    
    \begin{remark}
        An interesting feature of Theorem \ref{1.2} is that for any $c>1$ with $c\not\in\N$ the set $D_c=\bigcap_{i=1}^k \Delta_{2,c} (S_i)$ is thick. One can show that if $c\not\in\N$, unlike in Theorem \ref{c in N}, then the set $D_c$ may not be syndetic. Also when $c\in \N$ the syndetic set $D_c=\bigcap_{i=1}^k \Delta_{2,c} (S_i)$ in Theorem \ref{c in N} may not be thick. A detailed discussion regarding syndeticity and thickness of sets of the form $D_g=\bigcap_{i=1}^n \Delta_{2,g}(S_i)$ for $g:\N\to [1,\oo)$, that appear in Theorems \ref{3.9} and \ref{3.13} (which contain Theorems \ref{1.2} and \ref{c in N} as rather special corollaries), is presented in subsection \ref{Sec 3.3}. 
    \end{remark}

    As mentioned above, one of the major goals of this paper is to demonstrate that ergodic methods allow one to amplify and generalize Theorem \ref{Ruz} to the setup of countable cancellative amenable semigroups. Theorem \ref{1.3}, formulated below and proved in Section \ref{Sec 4}, generalizes Theorem \ref{Ruz} in more than one way. First, Theorem \ref{1.3} extends the setup from $\Z$ to a general countably infinite amenable group $G$. Second, for the sets $S_1,\dots S_k$, which appear in the formulation of Theorem \ref{1.3}, it is assumed $\ud_{\mathbf{F}_i}(S_i) >0$ for $i=1,\dots k$, where the F\o lner sequences $\mathbf{F}_1,\dots,\mathbf{F}_k$ are arbitrarily chosen. Finally, the tempered F\o lner sequence $\mathbf{G}$, which appears in the formulation of Theorem \ref{1.3}, leads even in the case $G=\Z$ to a nontrivial generalization of Theorem \ref{Ruz}. We will discuss Følner and tempered Følner sequences in detail in subsection \ref{Sec 4.2}.
    

\begin{Theorem}[ (Theorem \ref{general ruzsa}, subsection \ref{Sec 4.2})] \label{1.3}
    Let $G$ be a countably infinite amenable group with left Følner sequences $\mathbf{G},\mathbf{F}_1,\ldots,\mathbf{F}_k$, where $\mathbf{G}$ is left tempered. Let $S_1,\ldots,S_k \sseq G$ such that $\overline{d}_{\mathbf{F}_i}(S_i)> 0$ for each $i=1,\ldots,k$. Then there is $S \sseq G$ such that $d_\mathbf{G}(S) \geq \prod_{i=1}^k \overline{d}_{\mathbf{F}_i}(S_i)$ and $
                \Delta_1(S) \sseq D = \bigcap_{i=1}^k \Delta_2(\mathbf{F}_i,S_i)$. Furthermore, there are $m_1,\ldots,m_\ell \in G$, where
            \begin{equation}
                \ell \leq \prod_{i=1}^k 1/\overline{d}_{\mathbf{F}_i}(S_i),
            \end{equation}
        such that $\bigcup_{i=1}^\ell(m_iD) = \bigcup_{i=1}^\ell(Dm_i^{-1}) = G$. 
\end{Theorem}


    The structure of this paper is as follows. The goal of Section \ref{Sec 2} is to provide a short ergodic proof of Theorem \ref{Ruz} and to set up some ideas that will be further developed and amplified in subsequent sections. In Section 3, we obtain Theorems \ref{3.9} and \ref{3.13} dealing with the difference sets defined in (\ref{1.6}), which contain, respectively, Theorems \ref{1.2} and \ref{c in N} as instantaneous corollaries. In Section 4, we will prove Theorem \ref{1.3} and its more general version, Theorem \ref{general ruzsa for semigroups}, which pertains to cancellative amenable semigroups. Finally, in Section 5, we give examples of tempered Følner sequences in various semigroups such as $(\N,+)$, $(\N,\cdot)$, the Heisenberg group, and locally finite groups. 

    \textbf{Acknowledgments:} The authors thank John Griesmer and Saúl Rodríguez Martín for helpful communications.

\section{Ergodic proof of Ruzsa's Theorem} \label{Sec 2}
The proof of Theorem \ref{Ruz} we present in this section will utilize two facts from ergodic theory. The first of them is an extension of the classical Poincaré recurrence theorem.\footnote{Poincaré's recurrence theorem states that if $(X, \ca{B},\mu)$ is a probability space, $T$ a $\ca{B}$-measurable map $X\to X$ such that $\mu(T^{-1}A) = \mu(A)$ for all $A \in \ca{B}$, and $B \in \ca{B}$ with $\mu(B) > 0$, then there is $n \in \N$ such that $\mu(B \cap T^{-n}B) > 0$.}


\begin{Theorem}[ (\site{Theorem 1.2}{Ber85})] \label{Berg}
       Let $(X, \mathcal{B}, \mu, T)$ be a measure-preserving system\footnote{A \textit{measure-preserving system} is a quadruple $(X,\ca{B},\mu,T)$, where $(X,\ca{B},\mu)$ is a probability space and $T$ is a $\ca{B}$-measurable map $X\to X$ such that for any $A\in\ca{B}$, one has $\mu(T^{-1}A)=\mu(A)$. A measure-preserving system $(X,\ca{B},\mu,T)$ is \textit{invertible} if the map $T$ is invertible. }
       and let $B \in \mathcal{B}$ with $\mu(B) > 0$. Then there is a set $P \subseteq \N$ with $d(P) \geq \mu(B) $ such that for any $n_1,\ldots,n_r \in P$ we have $\mu \left( B \cap T^{-n_1}B \cap \cdots \cap T^{-n_r}B \right) > 0$. 
    \end{Theorem}

    \begin{remark}
    \label{2.2}
        We would like to stress a subtle point in the statement of Theorem \ref{Berg}, namely that the set $P$ has positive \textit{natural density}. This fact is a consequence of 
        the following variant of the classical pointwise ergodic theorem: for any invertible measure-preserving system $(X,\mathcal{B}, \mu, T)$ and bounded measurable function $f:X \rightarrow \R$, the limit 
        \begin{equation}\label{eq 2.1}
            \lim_{N \rightarrow \infty} \frac{1}{N} \sum_{n=1}^N f (T^{n}x)
        \end{equation}
            exists for almost every $x \in X$.
    \end{remark}



The second fact that we will need is a variant of Furstenberg's correspondence principle. 

\begin{Theorem}[ (cf. \site{Theorem 1.1}{Ber87})] \label{Furst}
    Suppose that $S \sseq \N$ with $\overline{d}(S) > 0$. Then there is an invertible measure-preserving system $(X,\mathcal{B},\mu,T)$ and a set $B \in \mathcal{B}$ with $\mu(B) = \overline{d}(S)$ such that
        \begin{equation}\label{Furst A}
            \overline{d}\left(S \cap (S-n_1) \cap \cdots \cap (S-{n_k}) \right)\geq\mu\left( B \cap T^{-n_1}B \cap \cdots \cap T^{-n_k}B \right) 
        \end{equation}
    for all $n_1,\ldots,n_k \in \Z$. 
\end{Theorem}


        \noindent \textit{Proof of Theorem \ref{Ruz}.}
            By Theorem \ref{Furst}, for each $i = 1,\ldots,k$ there exists an invertible measure-preserving system $(X_i, \mathcal{B}_i, \mu_i, T_i)$ and a set $B_i \in \mathcal{B}_i$ with  $\mu_i(B_i) = \overline{d}(S_i)$ such that $S_i$, $B_i$, and $T_i$ satisfy formula (\ref{Furst A}). Let $(X,\mathcal{B}, \mu,T)$ be the product system where $X = X_1 \times \cdots \times X_k$, $\mathcal{B}$ is the $\sigma$-algebra on $X$ generated by sets of the form $A_1 \times \cdots \times A_k$ with $A_i \in \mathcal{B}_i$ for each $i= 1,\ldots,k$, $\mu$ is the product measure on $\ca{B}$, and $T:X \rightarrow X$ is the product transformation defined by $T(x_1,\ldots,x_k) = (T_1(x_1),\ldots,T_k(x_k))$. \\
            \indent Let $B = B_1\times \cdots \times B_k$. 
            By Theorem \ref{Berg} there exists $S\subseteq \N$ such that $d(S)\geq \mu(B)$ and for any $n_1,\ldots,n_r \in S$ we have
                \begin{equation} \mu \left( B \cap T^{-n_1}B \cap \cdots \cap T^{-n_r}B \right) > 0.\end{equation}
            Note that
                \begin{equation} d(S) \geq \label{2.4}\mu(B)=\prod_{i=1}^k\mu_i(B_i)=\prod_{i=1}^k\ud(S_i).\end{equation}
            Let $D$ denote the set $\bigcap_{i=1}^k\Delta_2 (S_i)$. In order to show that $\Delta_1 (S)\subseteq D$, take any $n-m \in \Delta_1(S)$ where $n,m\in S$ and note that for all $i\in\set{1,...,k}$, one has
            \begin{align}
                0<\mu\leri{ T^{-n}B\cap T^{-m}B }&=\mu\leri{B\cap T^{-(n-m)}B}\\
                &\le \mu_i\leri{B_i\cap T_i^{-(n-m)}B_i}\le \ud(S_i\cap(S_i-(n-m))).
            \end{align}
            It follows that $n-m\in \Delta_2 (S_i)$ for all $i\in\set{1,...,k}$, which implies that $n-m\in D$. Thus, $\Delta_1 (S)\subseteq D$ as desired. This concludes the ergodic portion of our proof. \\
            \indent For the sake of completeness, we now prove the claim regarding the syndeticity of $D$
            in a way similar to 
            the proof of Theorem 2 in \cite{Ruz78} 
            (the only difference is that while the argument there is phrased in terms of so-called \textit{homogeneous systems}, we avoid them entirely).\\
            \indent One can observe that for any distinct $a_1,...,a_r\in\N$ with $r$ sufficiently large, two of the sets $S-a_1,...,S-a_r$ must have non-empty intersection since $d(S)>0$ (in fact one only needs $\ud(S)>0$). It follows that for some finite $F\subseteq S$, two of the sets $F-a_1,...,F-a_r$ have non-empty intersection. Indeed, if $s\in (S-a_i)\cap (S-a_j)$ for some distinct $i,j\in\set{1,...,r}$, one can take $F=\set{s+a_i,s+a_j}$.\\ 
            \indent Thus, there exists a set
            $M=\set{m_1,...,m_\ell}\subseteq\Z$ 
            of maximal cardinality $\ell$ such that for any finite $F\subseteq S$, the sets $F+m_1,...,F+m_\ell$ are pairwise disjoint. 
            Now let $m=\max_{1\le i\le \ell}|m_i|$ and note that for any $n\in\N$, the set $S(n)\defeq S\cap\set{1,...,n}$ satisfies
            \equat{
                \bigcup_{i=1}^\ell (S(n)+m_i)\subseteq [1-m,n+m].
            }
            Since the sets in the above union are disjoint, it follows that $\ell |S(n)|\le n+2m$. Dividing both sides by $n$ and taking the limit as $n\to\infty$ yields $\ell d(S)\le 1$, whence it follows from (\ref{2.4}) that
            $\ell\le 1/d(S)\le 1/\prod_{i=1}^k \ud(S_i)$. \\
            \indent We will now show that $\Z=\bigcup_{i=1}^\ell (D+m_i)$. To that end, let $n\in\Z$ be arbitrary. By maximality of $M$, there exists $i\in\set{1,...,\ell}$ and a finite $T\subseteq S$ with $(T+n)\cap (T+m_i)\neq\emptyset$. It follows that there exist $t,t'\in T$ for which $t+m_i=t'+n$. Thus, $n=t-t'+m_i\in \Delta_1 (S)+m_i\subseteq D+m_i$, and this concludes the proof.
            \qed

        

\section{Two general variants of Theorem \ref{Ruz} for subsets of the integers}

The goal of this section is to prove Theorems \ref{3.9} and \ref{3.13}, each of which can be seen as a general variant of Theorem \ref{Ruz} and which have, correspondingly, Theorems \ref{1.2} and \ref{c in N} as rather special corollaries. The proofs of Theorems \ref{3.9} and \ref{3.13} are presented in subsection \ref{Sec 3.2} and are based on an amplification of the ergodic ideas utilized in the previous section. Subsection \ref{Sec 3.2} is preambled by subsection \ref{ergodic sequences}, where we collect some basic facts which will be needed in the proof of Theorems \ref{3.9} and \ref{3.13}. 
Finally, in subsection \ref{Sec 3.3}, we collect some examples that illustrate the subtle features of Theorems \ref{3.9} and \ref{3.13} which pertain to the phenomena of thickness and syndeticity.



\subsection{Ergodic sequences} \label{ergodic sequences}

 As was mentioned in Remark \ref{2.2}, the proof of Theorem \ref{Berg}, and hence our proof of Theorem \ref{Ruz}, rely on the fact that the pointwise ergodic theorem holds along the sequence $1,2,3,\ldots$ of positive integers (see formula (\ref{eq 2.1})). There are actually many sequences $(a_n)_{n=1}^\infty$ in $\N$ with the property that for any invertible measure-preserving system $(X,\mathcal{B}, \mu,T)$ and bounded measurable function $f:X \rightarrow \R$, the limit
    \begin{equation} \label{eqn3.1}
            \lim_{N \rightarrow \infty} \frac{1}{N} \sum_{n=1}^N f (T^{a_n}x)
        \end{equation}
    exists for almost every $x \in X$. \\
\indent It turns out that by working with sequences satisfying formula (\ref{eqn3.1}) and by introducing some additional amplifications of the ideas presented in the proof of Theorem \ref{Ruz}, one can obtain new variants of Theorem \ref{Ruz}, such as Theorem \ref{1.2} \\
\indent Before formulating the next definition, we remind the reader that a measure-preserving transformation $T:X \rightarrow X$ on a probability measure space $(X,\mathcal{B}, \mu)$ is \textit{ergodic} if the only $T$-invariant sets (i.e.,
sets satisfying $1_A=1_{T^{-1}A}$ almost everywhere) are those with measure $0$ or $1$. 

\begin{definition} \label{ergodic defs}A sequence $(a_n)_{n=1}^\infty$ in $\N$ is 
\begin{itemize}[noitemsep,topsep=0pt]
    \item[i)] \textit{pointwise good} if for every measure-preserving system $(X, \ca{B}, \mu,T)$ and any function $f \in L^2$ the limit
    \begin{equation}
        \lim_{N \rightarrow \infty} \frac{1}{N} \sum_{n=1}^N f(T^{a_n}x)
    \end{equation}
    exists for almost every $x\in X$.
    \item[ii)] \textit{norm ergodic}\footnote{Although we will not use this result, we mention in passing that $(a_n)_{n=1}^\infty$ is norm ergodic if and only if for any Hilbert space $\mathcal{H}$ and unitary operator $U:\mathcal{H}\to \mathcal{H}$ which has no invariant elements in $\mathcal{H}$,
            \begin{equation} \label{7.1}
            \lim_{N \rightarrow \infty} \n{\frac{1}{N} \sum_{n=1}^N U^{a_n} f}_{\mathcal{H}} = 0.
        \end{equation}} 
    if for every ergodic invertible measure-preserving system $(X, \mathcal{B}, \mu, T)$ and 
    any function $f \in L^2$,
        \begin{equation}
        \label{norm ergodic}
            \lim_{N \rightarrow \infty} \n{\frac{1}{N} \sum_{n=1}^N T^{a_n} f - \int f}_2 = 0.
        \end{equation}
        
    \item[iii)] \textit{pointwise ergodic} if for every ergodic invertible measure-preserving system $(X, \ca{B}, \mu,T)$ and any function $f \in L^2$,
    \begin{equation}
        \lim_{N \rightarrow \infty} \frac{1}{N} \sum_{n=1}^N f (T^{a_n}x) = \int f
    \end{equation}
for almost every $x\in X$. 
\end{itemize}

    \end{definition}

    Note that if a sequence $(a_n)_{n=1}^\infty$ in $ \N$ is pointwise good and norm ergodic, then it is pointwise ergodic. Also, using ergodic decomposition, one can verify that if $(a_n)_{n=1}^\infty$ is pointwise ergodic, then it is simultaneously pointwise good and norm ergodic. Thus, when dealing with pointwise ergodic sequences, it is suitable to direct our attention to pointwise good sequences and norm ergodic sequences. \\
     \indent The following theorem provides a necessary and sufficient condition for a sequence to be norm ergodic.
    
    \begin{Theorem}[ (\site{ Corollary 3}{BB74})] \label{norm ergodic condition}
    A sequence $(a_n)_{n=1}^\infty$ in $\N$ is norm ergodic if and only if for all $\lambda\in \R\setminus\Z$,
        \begin{equation} \label{7.2}
            \lim_{N\rightarrow \infty} \frac{1}{N} \sum_{n=1}^N e^{2 \pi i \lambda a_n} =  0.
        \end{equation}
    \end{Theorem}



\begin{remark}
\label{iff norm ergodic}
    It follows from Theorem \ref{norm ergodic condition} along with the Weyl criterion (\cite[Chapter 1, Theorem 2.1]{KN74} and \cite[Chapter 5, Corollary 1.1]{KN74}) that a sequence $(a_n)_{n=1}^\oo$ in $\N$ is norm ergodic if and only if both of the following hold:
        \begin{itemize}[noitemsep, topsep=0pt]
        \item[i)] $(\lambda a_n)_{n=1}^\oo$ is uniformly distributed  mod 1 for any irrational $\lambda$.
        \item[ii)] $(a_n)_{n=1}^\oo$ is uniformly distributed in $\Z$, meaning that, for all $m \in \Z$,
            \begin{equation}
        \lim_{N \rightarrow \infty} \frac{|\{1 \leq n \leq N: a_n \equiv j \text{ mod } m\}|}{N} = 1/m \quad \text{ for } j = 1,\ldots,m. 
    \end{equation}
    \end{itemize} 
    \indent Equivalently, $(a_n)_{n=1}^\infty$ is norm ergodic if and only if formula (\ref{norm ergodic}) holds for two types of measure-preserving systems: irrational rotations on the unit circle, and transformations of the form $x \mapsto x+1$ on $\Z/k\Z$ for $k \in \N$.
\end{remark}

In what follows, we will be working with sequences of the form $([g(n)])_{n=1}^\infty$ for some natural classes of functions $g:\N \rightarrow \R$. The following theorem provides a useful sufficient condition for a sequence $([g(n)])_{n=1}^\infty$ to be norm ergodic. 

\begin{Theorem}
\label{3.4}
    Let $g:\N \rightarrow [1,\infty)$ be a function. If $(\lambda g(n))_{n=1}^\infty$ is uniformly distributed mod 1 for all irrational $\lambda$ and $(g(n)/m)_{n=1}^\infty$ is uniformly distributed mod 1 for all $m \geq 2$, then $([g(n)])_{n=1}^\infty$ is norm ergodic.
\end{Theorem}

\begin{Proof}
    Imitating the proof of \cite[Theorem 5.12]{BK09}, we will first prove that $([g(n)]\gamma)_{n=1}^\infty$ is uniformly distributed mod 1 when $\gamma$ is irrational. For all $(a,b) \in \Z^2 \setminus\{(0,0)\}$ the sequence $(ag(n)\gamma+bg(n))_{n=1}^\infty=((a\gamma+b)g(n))_{n=1}^\infty$ is uniformly distributed mod 1, so $(g(n)\gamma,g(n))_{n=1}^\infty$ is uniformly distributed mod 1 (\cite[Chapter 1, Theorem 6.3]{KN74}). Define $f(x,y):=e^{2\pi i h(x-\{y\}\gamma)}$ and observe that, for all $h \in \Z \setminus \{0\}$, we have
        \begin{align*}
        \lim_{N\to\oo}\frac{1}{N}\sum_{n=1}^N e^{2\pi i h[g(n)]\gamma}& =\lim_{N\to\oo}\frac{1}{N}\sum_{n=1}^N e^{2\pi i h (g(n)\gamma-\{g(n)\}\gamma)} =\lim_{N\to\oo}\frac{1}{N}\sum_{n=1}^Nf(g(n)\gamma,g(n)),
        \end{align*}
    where $\{\cdot\}$ denotes the fractional part. Since $f(x,y)$ is a Riemann-integrable periodic mod 1 function on $\R^2$,
    \[\lim_{N\to\oo}\frac{1}{N}\sum_{n=1}^N e^{2\pi i h[g(n)]\gamma}=\lim_{N\to\oo}\frac{1}{N}\sum_{n=1}^Nf(g(n)\gamma,g(n)) = \int_0^1\int_0^1f(x,y)dxdy=0.\]
    By the Weyl criterion (\cite[Chapter 1, Theorem 2.1]{KN74}) it therefore follows that $([g(n)]\gamma)_{n=1}^\infty$ is uniformly distributed mod 1. \\
\indent By \cite[Chapter 5, Theorem 1.4]{KN74}, since $(g(n)/m)_{n=1}^\infty$ is uniformly distributed mod 1 for all integers $m\geq 2$ the sequence $([g(n)])_{n=1}^\infty$ is uniformly distributed in $\Z$. It follows from Remark \ref{iff norm ergodic} that $([g(n)])_{n=1}^\infty$ is norm ergodic.
\end{Proof}
We will now describe a rather general class of functions which, with the help of Theorem \ref{3.4}, provide a wide variety of norm ergodic sequences of the form $([g(n)])_{n=1}^\infty$.\\
\indent Before stating the next theorem we first introduce the notion of a \textit{Hardy field}. Let $B$ be the set of all germs at infinity\footnote{A \textit{germ at infinity} is any equivalence class of real-valued functions in one real variable under the equivalence relation $f \sim g$ if and only if there is $t_0 > 0$ such that $f(t) = g(t)$ for all $ t \geq t_0$.} of continuous real functions defined on $[1,\infty)$. Note that $B$ forms a ring with respect to pointwise addition and multiplication. A \textit{Hardy field} is any subfield of $B$ which is closed under differentiation. We will denote the union of all Hardy fields by $\mathbf{U}$. A function $g \in \mathbf{U}$ is called \textit{subpolynomial} if there is some positive integer $k$ such that $g(x)/x^k \rightarrow 0$ as $x \rightarrow \infty$. 


\begin{Theorem}[ (\site{Theorem 1.3}{Bos94})]
\label{3.5}
    Let $g\in \mathbf{U}$ be a subpolynomial function. Then the following two conditions are equivalent:
    \begin{enumerate}[noitemsep, topsep=0pt]
        \item[i)] The sequence $(g(n))_{n=1}^\oo$ is uniformly distributed mod 1.
        \item[ii)] For every $p(x)\in \Q[x]$,
        \begin{equation} \label{bosh condition}\lim_{x\to\oo}\frac{g(x)-p(x)}{\log(x)}=\pm\oo.\end{equation}
    \end{enumerate}
\end{Theorem}

Now, with the help of Theorem \ref{3.5}, one can use Theorem \ref{3.4} to produce plenty of examples of norm ergodic sequences. For example, $[\log^t(n)] \text{ with } t>1$ is norm ergodic.
Also, if $g\in\mathbf{U}$ and $k\in\N$ are such that $g(x)/x^{k} \rightarrow \oo$ but $g(x)/x^{k+1} \rightarrow 0$ as $x \rightarrow \infty$, then $\lambda g$ satisfies condition (ii) of Theorem \ref{3.5} for all nonzero $\lambda \in \R$ (this follows from repeated applications of L'Hospital), and so $([g(n)])_{n=1}^\infty$ is norm ergodic by Theorem \ref{3.4}. Some examples obtained in this way are
    \[ [n^c]\text{ with } c>0, c \notin \N,\quad [n\log(n)],\quad \left[\frac{n^2}{\log(n)}\right], \quad \quad \text{for } n = 1,2,3,\ldots \]
\indent To describe another class of examples, suppose that $f(x) \in \R[x]$ is a polynomial with two coefficients other than the constant term, which are linearly independent over $\Q$. Then $\lambda f$ satisfies condition (ii) of Theorem \ref{3.5} for all nonzero $\lambda \in \R$, and so $([f(n)])_{n=1}^\infty$ is norm ergodic by Theorem \ref{3.4}. So, for instance, the sequence $([\sqrt{2}n^2+\sqrt{3}n])_{n=1}^\infty$ is norm ergodic. \\
\indent It was shown in \cite[Theorem B]{BGAM05} that many sequences guaranteed to be norm ergodic by \cite{Bos94} and Theorem \ref{3.4} are also pointwise good. For example, the following sequences are pointwise good:
\[n^3+n^2+n+1,\quad \left[ n\log(n)\right],\quad \left[ n^c \right] \text{ with } c>0, \quad \left[\frac{n^2}{\log(n)}\right], \quad \quad \text{for } n = 1,2,3,\ldots\]
\indent On the other hand, one can show that, while $([\log^t(n)])_{n=1}^\infty$ is norm ergodic for any $t>1$, it is not pointwise good (this follows from \cite[Theorem 2.16]{JW94}; see \cite[Example 2.18]{JW94} for a proof of the case $t=2$).
\\
\indent As was mentioned above, if a sequence $(a_n)_{n=1}^\infty$ in $\N$ is pointwise good and norm ergodic, then it is pointwise ergodic. On the base of this, we can obtain many sequences which are pointwise ergodic. For example, the following sequences are pointwise ergodic:
    \begin{equation}
    \label{PW EX}
        [\sqrt{2}n^2+\sqrt{3}n], \quad [n^c] \text{ with } c>0, c \notin \N, \quad [n\log(n)], \quad \left[ \frac{n^2}{\log(n)} \right], \quad \quad \text{for }  n = 1,2,3,\ldots
    \end{equation}  

\subsection{Proofs of Theorems \ref{3.9} and \ref{3.13}}\label{Sec 3.2}


The goal of this subsection is to present a proof of a general variant of Theorem \ref{Ruz}, namely Theorem \ref{3.9}. The proof of Theorem \ref{3.9} will rely on two results that we will presently formulate and prove.

\begin{Theorem}[ (cf. \site{Corollary 7.2}{BK09})] \label{3.6}
        Suppose that $(a_n)_{n=1}^\infty$ is a norm ergodic sequence in $\N$. If $(X,\mathcal{B}, \mu, T)$ is an invertible measure-preserving system and $A \in \mathcal{B}$ with $\mu(A) > 0$, then   
                \begin{equation} \label{14}
                    \lim_{N \rightarrow \infty} \frac{1}{N} \sum_{n=1}^N \mu(A \cap T^{a_n} A) \geq \mu(A)^2.
                \end{equation} 
    \end{Theorem}

    \begin{Proof}
        By using weak convergence, if $T$ is ergodic, it follows from the assumption that $(a_n)$ is a norm ergodic sequence that
        \[\lim_{N\to\oo}\frac{1}{N}\sum_{n=1}^N\mu(A\cap T^{a_n}A)=\lim_{N\to\oo}\frac{1}{N}\sum_{n=1}^N\int 1_A T^{-a_n} 1_A d\mu=\left(\int 1_A \right)^2=\mu(A)^2.\]
        If $T$ is not ergodic, then the inequality \ref{14} follows by utilizing the ergodic decomposition of $\mu$ (an example of such an application of the ergodic decomposition is delineated in \cite[Section 5, page 48]{BER2000}).
    \end{Proof}


     \begin{Theorem} \label{3.7}
        Let $g: \N \rightarrow [1,\infty)$ be a function such that the sequence $([g(n)])_{n=1}^\infty$ is  pointwise ergodic, let $(X, \mathcal{B}, \mu,T)$ be an invertible measure-preserving system, and let $B \in \mathcal{B}$ with $\mu(B)> 0$. Then there is a set $P \sseq \N$ with $d(P) \geq \mu(B)$ such that for any $n_1,\ldots, n_r \in P$ we have $\mu(B \cap T^{-[g(n_1)]}B \cap \cdots \cap T^{-[g(n_r)]}B) > 0$. 

     \end{Theorem}

     \begin{Proof} We will employ an argument similar to the one used in the proof of \cite[Theorem 1.2]{Ber85}, cited above as Theorem \ref{Berg}. \\
     \indent After deleting (if needed) a set of measure zero from each $T^{-n}B$, we can and will assume without loss of generality that for any $n_1,...,n_j\in\Z$, if $T^{-n_1}B\cap\cdots\cap T^{-n_j}B\neq\emptyset$, then $\mu(T^{-n_1}B\cap\cdots\cap T^{-n_j}B)>0$.\\
     \indent For all $N \in \N$ define $f_N:X \rightarrow \N$ by
            \[ f_N(x) = \frac{1}{N} \sum_{n=1}^N 1_B\left(T^{[g(n)]}(x)\right) \]
        and let $f = \limsup_{N \rightarrow \infty} f_N$. Note that $0 \leq f_N \leq 1$ for all $N$, so $f$ is well-defined and takes values in $[0,1]$. We see that
            \begin{align} 
                \int 1_B \cdot f_N &= \int 1_B \cdot \frac{1}{N} \sum_{n=1}^N 1_B \left(T^{[g(n)]}\right) = \frac{1}{N} \sum_{n=1}^N \int 1_B \cdot 1_B\left(T^{[g(n)]}\right) \\
                &= \frac{1}{N} \sum_{n=1}^N \int 1_{B} \cdot 1_{ T^{-[g(n)]}B} = \frac{1}{N} \sum_{n=1}^N \mu(B \cap T^{-[g(n)]} B), 
            \end{align}
        so by Theorem \ref{3.6},
            \begin{equation} \label{int}
                \lim_{N \rightarrow \infty} \int 1_B \cdot f_N = \lim_{N\rightarrow \infty}\frac{1}{N} \sum_{n=1}^N \mu(B \cap T^{-[g(n)]} B) \geq \mu(B)^2. 
            \end{equation}
        Noting that $1_B\cdot f = \limsup_{N\rightarrow \infty}(1_B \cdot f_N)$, Fatou's lemma gives us 
            \begin{equation} \label{Fatou}
               \int_B f = \int 1_B \cdot f = \int \limsup_{N \rightarrow \infty} (1_B \cdot f_N )\geq \limsup_{N \rightarrow \infty} \int 1_B \cdot f_N = \lim_{N \rightarrow \infty} \int(1_B \cdot f_N) \geq \mu(B)^2,
            \end{equation}
        thus
            \begin{equation} 
                \int_B f \geq \mu(B)^2 = \int_B \mu(B), 
            \end{equation}
        and so $\mu\left(B \cap \{f \geq \mu(B)\}\right) > 0$. 
        The sequence $([g(n)])_{n=1}^\infty$ is pointwise ergodic, hence pointwise good, so the limit
            \begin{equation} \label{lim}
                 \lim_{N \rightarrow \infty} \frac{1}{N}\sum_{n=1}^N 1_B\left(T^{[g(n)]}(x)\right) = \lim_{N \rightarrow \infty} f_N(x) 
            \end{equation}
        exists almost everywhere, so there must be some $x_0 \in B$ such that $f(x_0) \geq \mu(B)$ and $\lim_{N\rightarrow \infty} f_N(x_0)$ exists. Therefore we have 
            \begin{align} 
                \lim_{N \rightarrow \infty} \frac{1}{N} \sum_{n=1}^N 1_B(x_0) \cdot 1_B\left( T^{[g(n)]}(x_0)\right) &= 1_B(x_0) \cdot \lim_{N\rightarrow \infty} f_N(x_0) \\
                \label{3.16} &=  \lim_{N\rightarrow \infty} f_N(x_0) = f(x_0) \geq \mu(B).
            \end{align}
        Finally, putting $P = \{n \in \N: x_0 \in B \cap T^{-[g(n)]}B \}$ gives us what we want. Indeed, we have 
            \begin{align}
                \lim_{N \rightarrow \infty} 
                \frac{| P \cap \{1,\ldots,N\}|}{N} =\lim_{N \rightarrow \infty} \frac{1}{N} \sum_{n=1}^N 1_B(x_0) \cdot 1_B\left(T^{[g(n)]}(x_0)\right) \geq \mu(B)
            \end{align}
        by formula \ref{3.16}. Now suppose that $n_1,\ldots,n_m \in P$. Then for $i = 1,\ldots,m$ we have $x_0 \in B \cap T^{-[g(n_i)]}B$, so $B \cap T^{-[g(n_1)]}B \cap \cdots \cap T^{-[g(n_m)]}B \neq \emptyset$, thus $\mu ( B \cap T^{-[g(n_1)]}B \cap \cdots \cap T^{-[g(n_m)]}B) > 0$.
    \end{Proof}

    \begin{Corollary} \label{cor3.8}
        Let $g:\N \rightarrow [1,\infty)$ be a function so that the sequence $([g(n)])_{n=1}^\infty$ is pointwise ergodic and let $S_1,\dots,S_k\subseteq \N$ such that $\ud(S_i)>0$ for $1\leq i\leq k$. Then there exists a set $S \sseq \N$ with $d(S) \geq \prod_{i=1}^k \overline{d}(S_i)$ such that for all $n_1,\ldots,n_r \in S$ and for any $1 \leq i \leq k$ we have $\overline{d}(S_i \cap (S_i - [g(n_1)]) \cap \cdots \cap (S_i - [g(n_r)])) > 0$.
    \end{Corollary}

    \begin{Proof}
        By Theorem \ref{Furst}, for each $i = 1,\ldots,k$ there exists an invertible measure-preserving system $(X_i, \mathcal{B}_i, \mu_i, T_i)$ and a set $B_i \in \mathcal{B}_i$ with  $\mu_i(B_i) = \overline{d}(S_i)$ such that $S_i$, $B_i$, and $T_i$ satisfy formula (\ref{Furst A}). Let $(X,\mathcal{B}, \mu,T)$ be the product of the systems $(X_i,\mathcal{B}_i,\mu_i,T_i)$ for $1 \leq i \leq k$, where $B=B_1\times\dots\times B_k$ (see the proof of Theorem \ref{Ruz} for an explicit description of $(X,\mathcal{B}, \mu,T)$).  
        By Theorem \ref{3.7} there exists a set $S \sseq \N$ with $d(S) \geq \mu(B)$ such that for any $n_1,\ldots, n_r \in S$ we have $\mu(B \cap T^{-[g(n_1)]}B \cap \cdots \cap T^{-[g(n_r)]}B) > 0$. Then for any $1 \leq i \leq k$ and for all $n_1,\dots,n_r \in S$, we have 
        \begin{align}
           & \overline{d}(S_i \cap (S_i - [g(n_1)])\cap \dots\cap (S_i - [g(n_r)])) \geq \prod_{i=1}^k \overline{d}(S_i \cap (S_i - [g(n_1)])\cap \dots\cap (S_i - [g(n_r)]))\nonumber \\
             &\geq \prod_{i=1}^k \mu_i(B_i \cap T_i^{- [g(n_1)]}B_i \cap \dots\cap T_i^{- [g(n_r)]}B_i )\nonumber =\mu(B \cap T^{- [g(n_1)]}B \cap \dots\cap T^{- [g(n_r)]}B ) > 0.
        \end{align}
    \end{Proof}

    \indent We now embark on the proof of a generalization of Theorem \ref{Ruz} which has Theorem \ref{1.2} as an immediate corollary. First, given a function $g:\N \rightarrow [1,\infty)$, we define sets $\Delta_{1,g}(S)$ and $\Delta_{2,g}(S)$ which form a general version of the sets $\Delta_1(S)$ and $\Delta_2(S)$ defined in the introduction. 

    \begin{definition}
        For $S \subseteq \N$ and for $g:\N \rightarrow [1,\infty)$, let
            \begin{equation*} 
                \Delta_{1,g}(S) = \{n \in \N : S \cap (S-[g(n)]) \neq \emptyset\},\quad 
                \Delta_{2,g}(S) = \{n \in \N : \overline{d}(S \cap (S-[g(n)])) > 0 \}. 
            \end{equation*}
     \end{definition}
     \begin{Theorem}
     \label{3.9}
        Let $g:\N \rightarrow [1,\infty)$ and assume $([g(n)])_{n=1}^\infty$ is pointwise ergodic. If $S_1,\dots,S_k\subseteq \N$ satisfy $\ud(S_i)>0$ for each $i=1,\dots,k$, then there is $S\subseteq \N$ with $d(S)\geq \prod_{i=1}^k \ud(S_i)$ and $ \Delta_{1,g} (S)\subseteq D_g:=\bigcap_{i=1}^k \Delta_{2,g} (S_i)$.
        Furthermore, we have $\ud(\Delta_{1,g}(S))\geq \prod_{i=1}^k \ud(S_i)$. If, in addition, the function $g$ comes from a Hardy field with the property that for some $\ell\in\N\cup\{0\}$, $\lim_{x\to\oo}g^{(\ell)}(x)=\pm\oo$ and $\lim_{x\to\oo}g^{(\ell+1)}(x)=0$, then the set $D_g$ is thick. 
    \end{Theorem}

     \begin{Proof} By Theorem \ref{Ruz} there exists $S\subseteq \N$ with $d(S)\geq \prod_{i=1}^k \ud(S_i)$ and $\Delta_1(S)\subseteq \bigcap_{i=1}^k \Delta_2 (S_i)$. If $n\in \Delta_{1,g} (S)$, then $[ g(n)]\in \Delta_1(S)\subseteq\bigcap_{i=1}^k\Delta_2(S_i)$,which implies $n \in \bigcap_{i=1}^k \Delta_{2,g}(S_i)$. This shows that $\Delta_{1,g}(S) \sseq D_g$. \\
    \indent Next, by Corollary \ref{cor3.8}, there is a set $P \sseq \N$ with $d(P) \geq d(S)\geq  \prod_{i=1}^k \overline{d}(S_i)$ such that, for all $n \in P$ we have $\overline{d}(S \cap (S - [g(n)])) > 0$. It follows that $P \sseq \Delta_{1,g}(S)$, hence $\overline{d}(\Delta_{1,g}(S)) \geq d(P) \geq \prod_{i=1}^k \overline{d}(S_i)$.\\
    \indent Now assume in addition that the function $g$ comes from a Hardy field with the property that for some $\ell\in\N\cup\{0\}$, $\lim_{x\to\oo}g^{(\ell)}(x)=\pm\oo$ and $\lim_{x\to\oo}g^{(\ell+1)}(x)=0$. Then \cite[Theorem A]{BMR20} implies, via Furstenberg’s correspondence principle (Theorem \ref{Furst}), that the set $D_g$ is thick.\end{Proof}


    \indent Theorem \ref{3.9} specializes to Theorem \ref{1.2} when $g(n) = n^c$ for $c > 1, c \notin \N$.

    We conclude this subsection by formulating a variant of Theorem \ref{3.9} in which the condition of pointwise ergodicity of $([g(n)])$ is relaxed, which in turn extends the range of applicability of the result. But first, we need to introduce yet another definition.

    \begin{definition} \label{seq of recurrence} Let $(a_n)_{n=1}^\infty$ be a sequence in $\N$.
    \begin{itemize}
    \vspace{-0.5cm}
    \item[i)] $(a_n)$ is an \textit{averaging sequence of recurrence} if for every invertible measure-preserving system $(X,\mathcal{B},\mu,T)$ and any set $A \in \mathcal{B}$ with $\mu(A)>0$,
        \begin{equation}
            \lim_{N \rightarrow \infty} \frac{1}{N}\sum_{n=1}^N \mu(A \cap T^{a_n}A) >0.
        \end{equation} 
    \item[ii)]
    $(a_n)$ is a \textit{uniform averaging sequence of recurrence} if for every invertible measure-preserving system $(X,\mathcal{B},\mu,T)$ and any set $A \in \mathcal{B}$ with $\mu(A)>0$, 
                \begin{equation}
                    \lim_{N-M \rightarrow \infty} \frac{1}{N-M}\sum_{n=M+1}^N \mu(A \cap T^{a_n}A) >0.
                \end{equation} 
    \end{itemize}

Note that, by Theorem \ref{3.6}, any sequence which is norm ergodic is also an averaging sequence of recurrence. For example, $([n^c])_{n=1}^\infty$ for $c>0$ is an averaging sequence of recurrence.
On the other hand, one can show that when $c>0$, $c \notin \N$, the sequence $([n^c])_{n=1}^\infty$ is not a uniform averaging sequence of recurrence.
A good supply of examples of uniform averaging sequences of recurrence is provided by polynomials. For example, if $q(n) \in \Q[n]$ is such that $q(\Z) \sseq \Z$ and is \textit{intersective}, that is, if $\{q(n): n \in \Z\} \cap k\Z$
is non-empty for each $k \in \N$, then $(q(n))_{n=1}^\infty$ is a uniform averaging sequence of recurrence (see \cite[Theorem 1.2]{BKS}). Also, if $p(n) \in \R[n]$ has at least two coefficients other than the constant term which are linearly independent over $\Q$, then $([p(n)])_{n=1}^\infty$ is a uniform averaging sequence of recurrence (see \cite[Example 1.3]{BKS}). 
    \end{definition}

    %
    %


    \begin{Theorem} \label{3.13}
        Let $g:\N \rightarrow [1,\infty)$ and assume that $([g(n)])_{n=1}^\infty$ is an averaging sequence of recurrence and is pointwise good. If $S_1,\dots,S_k\subseteq \N$ satisfy $\ud(S_i)>0$ for each $i=1,\dots,k$, then there is $S\subseteq \N$ with $d(S)> 0$ and $ \Delta_{1,g} (S)\subseteq D_g:=\bigcap_{i=1}^k \Delta_{2,g} (S_i)$.
        Furthermore, we have that $\ud(\Delta_{1,g}(S))  \ge d(S)>0$. If, in addition, the sequence $([g(n)])_{n=1}^\infty$ is a uniform averaging sequence of recurrence, then the set $D_g$ is syndetic. 
    \end{Theorem}

    \begin{Proof} 
    Admittedly, the proof of Theorem \ref{3.13} has similarities to, and uses ideas from, the proofs of Theorem \ref{3.7}, Theorem \ref{3.9}, and Corollary \ref{cor3.8}. However, the proof is not just a ``mechanical" composition of the aforementioned proofs. In particular, the proof below has novel elements when it comes to showing that under the assumptions of Theorem \ref{3.13}, the set $D_g$ is syndetic.\\
    \indent By Theorem \ref{Furst}, for every $i=1,...,k$, there exists an invertible measure-preserving system $(X_i,\ca{B}_i,\mu_i,T_i)$ and a set $B_i\in\ca{B}_i$ with $\mu_i(B_i)=\ud(S_i)$ satisfying formula (\ref{Furst A}). Let $(X,\ca{B},\mu,T)$ be the product of the systems $(X_i,\ca{B}_i,\mu_i,T_i)$, $i=1,...,k$ and let $B=B_1\times\cdots\times B_k$.\\
    \indent We will utilize a method similar to the one used in the proof of Theorem \ref{3.7} to obtain the set $S$. Consider the family of sets $T^{-n}B$ with $n\in\Z$. After deleting (if needed) a set of measure zero from each $T^{-n}B$, we can and will assume without loss of generality that for any $n_1,...,n_j\in\Z$, if $T^{-n_1}B\cap\cdots\cap T^{-n_j}B\neq\emptyset$, then $\mu(T^{-n_1}B\cap\cdots\cap T^{n_j}B)>0$. Now for all $N\in\N$, define $f_N=\frac{1}{N}\sum_{n=1}^N 1_{T^{-[g(n)]}B}$ and put $f=\limsup_{N\to\infty}f_N$. Since $([g(n)])_{n=1}^\infty$ is an averaging sequence of recurrence, it follows that
    \begin{align}
        \lim_{N\to\infty}\int 1_B\cdot f_N=\lim_{N\to\infty}\frac{1}{N}\sum_{n=1}^N\int 1_B\cdot 1_{T^{-[g(N)]}B}=\lim_{N\to\infty}\frac{1}{N}\sum_{n=1}^N\mu(B\cap T^{-[g(n)]}B)>0.
    \end{align}
    Since $1_B\cdot f=\limsup_{N\to\infty} (1_B\cdot f_N)$, Fatou's lemma implies 
    \begin{align}
        \int_B f=\int 1_B\cdot f\ge\limsup_{N\to\infty}\int1_B\cdot f_N=\lim_{N\to\infty}\int1_B\cdot f_N>0.
    \end{align}
    Thus, $\mu(B\cap \set{f>0})>0$. Since $([g(n)])_{n=1}^\infty$ is pointwise good, the limit 
    \begin{align}
    \lim_{N\to\infty}f_N(x)=\lim_{N\to\infty}\frac{1}{N}\sum_{n=1}^N 1_{T^{-[g(n)]}B}(x)
    \end{align}
    exists for almost every $x\in X$. In particular, there exists $x_0\in B$ for which we have $f(x_0)=\lim_{N\to\infty}f_N(x_0)>0$. Now let $S=\set{n\in\N:x_0\in B\cap T^{-[g(n)]}B}$ and observe that (since $x_0 \in B$)
    \begin{align}
        d(S)=\lim_{N\to\infty}\frac{\abs{S\cap\set{1,...,N}}}{N}=
        \lim_{N\to\infty}\frac{1}{N}\sum_{n=1}^N1_{B\cap T^{-[g(n)]}B}(x_0) =\lim_{N\rightarrow \infty} f_N(x_0) = f(x_0)>0.
    \end{align}
    Moreover, for any $n_1,...,n_j\in S$, we have $x_0\in B\cap T^{-[g(n_1)]}B\cap\cdots\cap T^{-[g(n_j)]}B$, which implies $\mu(B\cap T^{-[g(n_1)]}B\cap\cdots\cap T^{-[g(n_j)]}B)>0$ in light of our earlier assumptions. Therefore, by invoking Furstenberg's Correspondence Principle (Theorem \ref{Furst}), we have that for any $\ell=1,...,k$,
    \equat{
        &\ud(S_\ell\cap (S_\ell-[g(n_1)])\cap\cdots\cap (S_\ell-[g(n_j)])\ge
        \prod_{i=1}^k \ud(S_i\cap (S_i-[g(n_1)])\cap\cdots\cap (S_i-[g(n_j)])\\
        &\geq \prod_{i=1}^k \mu_i(B_i\cap T^{-[g(n_1)]}
        B_i\cap\cdots\cap T^{-[g(n_j)]}B_i) =
        \mu(B\cap T^{-[g(n_1)]}
        B\cap\cdots\cap T^{-[g(n_j)]}B)>0.
        }
    In particular, for any $n\in S$ and $\ell=1,...,k$, one has $\ud(S_\ell\cap (S_\ell-[g(n)]))>0$. It follows that $\Delta_{1,g}(S)\subseteq D_g$. \\
    \indent Next, by Corollary \ref{cor3.8}, there is a set $P \sseq \N$ with $d(P) \geq d(S) >0$ such that, for all $n \in P$ we have $\overline{d}(S \cap (S - [g(n)])) > 0$. It follows that $P \sseq \Delta_{1,g}(S)$, hence $\overline{d}(\Delta_{1,g}(S)) \geq d(P) \ge d(S)> 0$. \\
    \indent Lastly, suppose that $([g(n)])$ is a uniform averaging sequence of recurrence. By Furstenberg's correspondence principle, there is an invertible measure-preserving system $(Y,\mathcal{A}, \nu,U)$ and a set $A \in \mathcal{A}$ with $\nu(A) = d(S)$ such that, for all $n_1,\ldots,n_r \in \N$, 
        \[ \overline{d}(S \cap (S - n_1) \cap \cdots \cap (S-n_r)) \geq \nu(A \cap U^{-n_1}A \cap \cdots \cap U^{-n_r}A). \]
    It follows that $1_{\Delta_{2,g}(S)}(n) \geq \overline{d}(S \cap (S-[g(n)])) \geq \nu(A \cap U^{-[g(n)]}A)$ for all $n \in \N$, and so
        \begin{align*}
            \liminf_{N-M \rightarrow \infty}& \frac{|\Delta_{1,g}(S) \cap \{M+1,\ldots,N\}|}{N-M} = \liminf_{N-M \rightarrow \infty} \frac{1}{N-M} \sum_{n=M+1}^N 1_{\Delta_{1,g}(S)} (n)   \\
            &\geq \liminf_{N-M \rightarrow \infty} \frac{1}{N-M} \sum_{n=M+1}^N 1_{\Delta_{2,g}(S)} (n) \geq \lim_{N-M \rightarrow \infty} \frac{1}{N-M} \sum_{n=M+1}^N \nu(A \cap U^{-[g(n)]}A) > 0.
        \end{align*}
    This shows that $\Delta_{1,g}(S)$ is syndetic, hence $D_g$ is as well since $\Delta_{1,g}(S) \sseq D_g$. 
    \end{Proof}


\subsection{Some subtle points about sets of the form $D_g$ in Theorems \ref{3.9} and \ref{3.13}}\label{Sec 3.3}

    An interesting feature of Theorems \ref{3.9} and \ref{3.13} is that they provide sufficient conditions for sets of the form $D_g$ to be thick or syndetic. 
    In particular, if the function $g$ comes from a Hardy field with the property that for some $\ell\in\N\cup\{0\}$, $\lim_{x\to\oo}g^{(\ell)}(x)=\pm\oo$ and $\lim_{x\to\oo}g^{(\ell+1)}(x)=0$, then the set $D_g$ is always thick, as mentioned in Theorem \ref{3.9}. Such is the case if $g$ is one of the functions $n\log(n)$, $\frac{n^2}{\log(n)}$, or $n^c$ for $c>0$ with $c\not\in \N$, which are listed in (\ref{PW EX}). 
    On the other hand, if $g$ is a Hardy field function such that $([g(n)])_{n=1}^\infty$ is pointwise ergodic and $g$ satisfies $\lim_{x\to\oo}g^{(\ell)}(x)\to K$ for some finite constant $K \neq 0$, the set $D_g$ in Theorem \ref{3.9} may be syndetic.
    For example, if $g(x)=a_nx^n+a_{n-1}x^{n-1}+\cdots+a_0\in\R[x]$ has at least two rationally independent coefficients other than the constant term,
    then $D_{g}$ is always syndetic (see Theorem \ref{3.13}).
    An example of a function $g:\N\to [1,\oo)$ (which satisfies $\lim_{x\to\oo}g'(x)\to 2$) and a set $S \sseq \N$ for which the set $D_g=\Delta_{2,g}(S)$ from Theorem \ref{3.9} need not be thick nor syndetic is provided below by Example \ref{Ex 3.13}.


\begin{example}\label{Ex 3.13}
    If $g=2n+2\sqrt{n}$, then $([g(n)])_{n=1}^\oo$ is pointwise ergodic and the set $D_g=\Delta_{2,g}(4\N)$ is neither thick nor syndetic.
\end{example}
\begin{Proof}
    It follows from \cite[Theorem 3.4]{BGAM05} that $([g(n)])_{n=1}^\oo$ is pointwise ergodic. We have
    \begin{align*}
        D_g=\{n\in\N:\ud(4\N\cap(4\N-[g(n)]))>0\} =\{n:[g(n)]\in 4\N\}
    \end{align*}
        Note that for any real number $t$ one has that $[t]\in 4\N$ if and only if $\{\frac{1}{4}t\}\in [0,\frac{1}{4})$. So,
    \begin{align}
        D_g=\left\{n:\left\{\frac{1}{4}g(n)\right\}\in \left[0,\frac{1}{4}\right)\right\} =\left\{n:\left\{\frac{1}{2}(n+\sqrt{n})\right\}\in \left[0,\frac{1}{4}\right)\right\}.
    \end{align}
    To see that this set is not thick, observe that if $n\in D_g$ and $\epsilon=\frac{1}{2}(\sqrt{n+1}-\sqrt{n})<\frac{1}{4}$, then 
    \begin{align}
        \left\{\frac{1}{2}((n+1)+\sqrt{n+1})\right\}=\left\{\frac{1}{2}+\frac{1}{2}(\sqrt{n+1}-\sqrt{n})+\frac{1}{2}(n+\sqrt{n})\right\}\in \left(\frac{1}{2},\frac{3}{4}+\epsilon\right).
    \end{align}
    So we have shown that if $n\in D_g$, then $(n+1)\not\in D_g$, which implies $D_g$ is not thick.\\
    \indent Now we will show that $D_g$ is not syndetic. Fix $0<\epsilon<\frac{1}{8}$, then there exist arbitrarily large $n\in \N$ such that $\left\{\frac{1}{2}(n+\sqrt{n})\right\}\in \left(\frac{1}{4},\frac{1}{4}+\epsilon\right)$. Now observe for any $k\in \N$ and for $n\in \N$ large enough such that $\left\{\frac{1}{2}(n+\sqrt{n})\right\}\in \left(\frac{1}{4},\frac{1}{4}+\epsilon\right)$ we have for all $0\leq i\leq k$ that $\left\{\frac{1}{2}(n+i+\sqrt{n+i})\right\}\in \left(\frac{1}{4},\frac{1}{4}+2\epsilon\right)\cup \left(\frac{3}{4},\frac{3}{4}+2\epsilon\right)$. This means for any $k\in \N$ we can find a set $\{n,n+1,\dots,n+k\}\subset A^c$, hence $A$ is not syndetic.
\end{Proof}

The following example demonstrates that the set $D_g$ in Theorem \ref{3.9}, while guaranteed to be thick, need not be syndetic.

 \begin{example}\label{Ex 3.14}
     Let $D_{3/2} = \Delta_{2,{3/2}}(2\N) = \{n \in \N: \overline{d}(2\N \cap (2\N - [n^{3/2}])) > 0\}$. Then $D_{3/2}$ is not syndetic. 
 \end{example}

\begin{Proof} Note that $D_{3/2} = \{n \in \N: [n^{3/2}] \in 2\N\} = \{n \in \N: \{ \frac{1}{2} n^{3/2}\} \in [0,1/2)\}$. In order to prove $D_{3/2}$ is not syndetic, let $N\in\N$ and observe that it suffices to show that there is $M \in \N$ for which $M, M+1,\ldots,M+N \notin D_{3/2}$, which is equivalent to showing that $\{ \frac{1}{2}(M+k)^{1/2} \} \in [1/2,1)$ for $k = 0,1,\ldots,N$. To this end, we first observe the following:
\begin{itemize}
    \item[(a)] $\left( \frac{1}{2}(x+1)^{3/2} - \frac{1}{2}x^{3/2}\right) - \frac{3}{4}\sqrt{x} \rightarrow 0$ as $x \rightarrow \infty$,
    \item[(b)] The sequence $\leri{\set{\frac{1}{2}n^{3/2}},\set{\frac{3}{4}\sqrt{n}}}$, $n\in\N$ is dense in $[0,1]^2$ (this sequence is actually uniformly distributed in $[0,1]^2$ -- this follows from Féjer's Theorem \cite[Corollary 2.1, page 14]{KN74} and the Weyl criterion \cite[Theorem 6.2, page 48]{KN74}),
    \item[(c)] $\sqrt{x+1} - \sqrt{x} \rightarrow 0$ as $x \rightarrow \infty$ which, along with observation (b), implies that  
    for any interval $(a,b)\subseteq [0,1]$, there exist arbitrarily large $m\in\N$ for which $\set{\frac{3}{4}\sqrt{m+k}}\in (a,b)$ for $k=0,1,...,N$.
\end{itemize}
Now, fix $\alpha<\beta$ such that $[\alpha,\beta] \sseq (0,\frac{1}{16})$, and choose $\delta > 0$ so that $(\frac{\alpha}{N} - \delta, \beta + \delta) \sseq (0,\frac{1}{16})$. In light of the above observations, we can choose $M \in \N$ satisfying
\begin{gather}\label{Conditions on M}
    \begin{split}\text{(i)}\,\,\, \left\{\frac{1}{2}M^{3/2}\right\} \in \leri{\frac{5}{8}, \frac{7}{8}},\quad \text{(ii)}\,\,\, \left\{\frac{3}{4}\sqrt{M+j}\right\} \in \leri{\frac{\alpha}{N}, \frac{\beta}{N}}\text{ for }j=0,1,...,N,\text{ and }\\
    \text{(iii)}\,\,\,\left( \frac{1}{2}(M+j)^{3/2} - \frac{1}{2}(M+j-1)^{3/2}\right) - \frac{3}{4}\sqrt{M+j-1}<\frac{\delta}{N}\text{ for }j = 1,\ldots,N.\end{split}
\end{gather}
Let $1\le k\le N$ and recall that it suffices to show $\set{\frac{1}{2}(M+k)^{3/2}}\in [1/2,1)$. Observe that
    \begin{align}
        &\left| \left(\frac{1}{2}(M+k)^{3/2} - \frac{1}{2}M^{3/2}\right)  - \sum_{i=1}^k \frac{3}{4}\sqrt{M+i-1} \right| \tag*{}\\
        &=  \left| \sum_{i=1}^k \left(\frac{1}{2}(M+i)^{3/2} - \frac{1}{2}(M+i-1)^{3/2}\right) - \sum_{i=1}^k \frac{3}{4}\sqrt{M+i-1} \right| \tag*{}\\
        &\leq \sum_{i=1}^k \left| \left(\frac{1}{2}(M+i)^{3/2} - \frac{1}{2}(M+i-1)^{3/2}\right) - \frac{3}{4}\sqrt{M+i-1} \right| < k \cdot \frac{\delta}{N} \leq N \cdot \frac{\delta}{N} = \delta. \label{3.14 eqn 1}
    \end{align}
By (ii) in formula (\ref{Conditions on M}),
it follows that $\left\{ \sum_{i=1}^k \frac{3}{4}\sqrt{M+i-1}\right\} \in (\frac{\alpha}{N},\beta)$, and so by formula (\ref{3.14 eqn 1}), we get that 
\begin{align}\label{frac_parts_are_close}
\set{\frac{1}{2}(M+k)^{3/2} - \frac{1}{2}M^{3/2}}=\set{\set{\frac{1}{2}(M+k)^{3/2}}-
\set{\frac{1}{2}M^{3/2}}}
\in \leri{\frac{\alpha}{N} - \delta, \beta + \delta} \sseq \leri{0,\frac{1}{16}}.\end{align}
Since fractional parts lie in $[0,1)$, it follows from formula (\ref{frac_parts_are_close}) that either 
\begin{gather*}
    \text{(1)}\,\,\, \set{\frac{1}{2}(M+k)^{3/2}} - \set{\frac{1}{2}M^{3/2}} \in \leri{-1,-\frac{15}{16}},\text{ or } \text{(2)}\,\,\, \set{\frac{1}{2}(M+k)^{3/2}} - \set{\frac{1}{2}M^{3/2}}\in \leri{0,\frac{1}{16}}.
\end{gather*} 
By (i) from formula (\ref{Conditions on M}), $ \{\frac{1}{2}(M+k)^{3/2}\} - \{\frac{1}{2}M^{3/2}\} > 0-\frac{7}{8} > -\frac{15}{16}$, so case (1) is not possible. Hence, case (2) must hold and by (i) from formula (\ref{Conditions on M}), so it follows that $\{\frac{1}{2}(M+k)^{3/2}\} \in (\frac{10}{16}, \frac{15}{16}) \sseq [\frac{1}{2}, 1)$ as desired. \end{Proof}

\section{Generalizing Theorem \ref{Ruz} to the Amenable Setup}\label{Sec 4}

The goal of this section is to develop some far-reaching extensions of Theorem \ref{Ruz}. In Section \ref{Sec 4.1} we obtain an amplification of Theorem \ref{Ruz} which applies to general families of ``large sets" $S_1,\dots, S_k\subset \Z$. In Section \ref{Sec 4.2} we generalize Theorem \ref{Ruz} to the setup of countable amenable groups via Theorem \ref{general ruzsa}. Finally, in Section \ref{Sec 4.3} we generalize Theorem \ref{general ruzsa} to countable amenable cancellative semigroups.

\subsection{Amplifications of Theorem \ref{Ruz} in $\Z$}\label{Sec 4.1}



The ergodic techniques that were used in the previous sections admit further amplifications, which will allow us to establish a rather general form of Theorem \ref{Ruz} for countably infinite amenable groups and cancellative semigroups. Actually, as we will see in Theorem \ref{4.1} below, the point of view based on the notion of amenability allows for interesting extensions of Theorem \ref{Ruz} even for subsets of $\Z$. Indeed, a careful examination of the proof of Theorem \ref{Ruz} in Section 2 reveals that the assumption that the sets $S_1,\ldots,S_k$ have positive upper density can be significantly weakened.\\
\indent Before formulating Theorem \ref{4.1}, we need to introduce some additional notation. Recall that a Følner sequence in $\Z$ is a sequence $(F_N)_{N=1}^\oo$ of finite non-empty subsets of $\Z$ satisfying
            \begin{equation}
                \lim_{N\rightarrow\infty} \frac{|F_N \cap (F_N + n)|}{|F_N|} = 1
            \end{equation}
for all $n \in \Z$.\\
\indent Let $\mathbf{F} = (F_N)_{N=1}^\oo$ be a Følner sequence in $\Z$. Given a set $S \sseq \Z$, define
            \begin{equation}
                \label{eq 4.2}\overline{d}_\mathbf{F}(S) = \limsup_{N\rightarrow \infty}\frac{|S \cap F_N|}{|F_N|}.
            \end{equation}
        If the
limit in formula (\ref{eq 4.2}) exists, then its value is called the \textit{density} of $S$ along $\mathbf{F}$ and is denoted by $d_{\mathbf{F}}(S)$. Additionally, let
        \begin{equation}
            \Delta_1(S) = \{n \in \Z: S \cap (S-n) \}\quad\text{and}\quad \Delta_2(\mathbf{F},S) = \{n \in \Z: \overline{d}_\mathbf{F}(S \cap (S-n)) > 0\}. 
        \end{equation}
        It is worth mentioning that it is often natural to consider subsets and F\o lner sequences in $\N$ instead of $\Z$; the above definitions admit trivial modifications which are applicable to this case. \\
    \indent Note that $\Delta_2(\mathbf{F},S)$ coincides with $\Delta_2(S)$ and $\ud_{\mathbf{F}}(S)$ coincides with $\ud(S)$ as defined in Section 1 when $\mathbf{F}$ is the standard Følner sequence $\{1,\ldots,N\}$, $N = 1,2,3,\ldots$ \\
    \indent We are now ready to formulate an amplified version of Theorem \ref{Ruz}.


    \begin{Theorem}
    \label{4.1}
        Let $\mathbf{F}_1,\ldots,\mathbf{F}_k$ be Følner sequences in $\Z$. If $S_1,\ldots,S_k \sseq \Z$ satisfy $\overline{d}_{\mathbf{F}_i}(S_i)> 0$ for each $i=1,\ldots,k$, then there exists $S \sseq \Z$ such that $d(S) \geq \prod_{i=1}^k \overline{d}_{\mathbf{F}_i}(S_i)$ and $
                \Delta_1(S) \sseq D := \bigcap_{i=1}^k \Delta_2(\mathbf{F}_i,S_i)$. Furthermore, there are $m_1,\ldots,m_\ell \in \Z$, where
            \begin{equation}
                \ell \leq \prod_{i=1}^k 1/\overline{d}_{\mathbf{F}_i}(S_i),
            \end{equation}
        and $\bigcup_{i=1}^\ell(D+m_i) = \Z$. 

    \end{Theorem}

One can also obtain an amplification of Theorem \ref{3.9} similar to Theorem \ref{4.1}. For any function $g:\N \rightarrow [1,\infty)$ and any Følner sequence $\mathbf{F}$ in $\N$, let
    \[ \Delta_{1,g}(S) = \{n \in \N: S \cap (S-[g(n)]) \neq \emptyset\}, \]
    \[\Delta_{2,g}(\mathbf{F},S) = \{n \in \N: \overline{d}_{\mathbf{F}}(S \cap (S - [g(n)]))>0\}. \] 

\begin{Theorem} \label{4.2}
    Let $g: \N \rightarrow [1,\infty)$ and assume $([g(n)])_{n=1}^\infty$ is pointwise ergodic. Let $\mathbf{F}_1,\ldots,\mathbf{F}_k$ be Følner sequences in $\N$. If $S_1,\ldots,S_k \sseq \N$ satisfy $\overline{d}_{\mathbf{F}_i}(S_i) > 0$ for each $i = 1,\ldots,k$, then there exists $S \sseq \N$ such that $d(S) \geq \prod_{i=1}^k \overline{d}_{\mathbf{F}_i}(S_i)$ and $\Delta_{1,g}(S) \sseq D_g := \bigcap_{i=1}^k \Delta_{2,g}(\mathbf{F}_i,S_i)$. Furthermore, we have that $\overline{d}(\Delta_{1,g}(S)) \geq \prod_{i=1}^k \overline{d}_{\mathbf{F}_i}(S_i)$. 
\end{Theorem}

We will not prove Theorem \ref{4.1} or Theorem \ref{4.2} here, as their proofs are practically the same as those of Theorems \ref{Ruz} and \ref{3.9}, respectively, with one modification which involves the use of a more general version of Theorem \ref{Furst} in which $\overline{d}$ is replaced by $\overline{d}_\mathbf{F}$ (this result is a special case of Theorem \ref{general Furst} which is formulated in the next subsection). Theorem \ref{4.1} is also a special case of Theorem \ref{general ruzsa}, which will be proved in the next subsection.  \\
\indent Note that in Theorem \ref{4.1}, besides the F\o lner sequences $\mathbf{F}_1,\ldots,\mathbf{F}_k$, one more F\o lner sequence is implicitly present. Namely, $d(S)$ actually means $d_\mathbf{F}(S)$, where $\mathbf{F}$ is the standard Følner sequence $\{1,\ldots,N\}$, $N = 1,2,3,\ldots$ The role of the F\o lner sequence $\mathbf{F}$ is different from the roles of the sequences $\mathbf{F}_1,\ldots,\mathbf{F}_k$. While $\mathbf{F}_1,\ldots,\mathbf{F}_k$ are needed just to define appropriate notions of largeness for $S_1,\ldots,S_k$, the role of the sequence $\mathbf{F}$ is to guarantee that the pertinent ergodic averages converge pointwise (which is instrumental in the proof of Theorem \ref{Berg}, formulated in Section \ref{Sec 2}). These ideas will be clarified in the next subsection where we extend Theorem \ref{4.1} to the general amenable setup. 

    \subsection{Extending Theorem \ref{Ruz} to amenable groups}\label{Sec 4.2}

    As was mentioned in the introduction, a natural framework for generalizing Theorem \ref{Ruz} is that of amenable semigroups. In this subsection, we first introduce some definitions pertaining to the notion of amenability and establish an extension of Theorem \ref{Ruz} for amenable groups. While being of interest on its own, this result will serve as a tool for extending Theorem \ref{Ruz} to amenable semigroups, which will be done in subsection \ref{Sec 4.3}.\\
    \indent A (discrete) semigroup $(G,\cdot)$ is said to be \textit{left amenable} if there exists a left invariant mean on $\ell_\infty(G)$, meaning a positive, linear functional 
    $\lambda\col\ell_\infty(G)\to\R$ satisfying
    \begin{enumerate}
        \item[i)] $\lambda(1_G)=1$, and
        \item[ii)] $\lambda({}_x\phi)=\lambda(\phi)$ for every $\phi\in\ell_\infty(G)$ and $x\in G$, where ${}_x\phi(t)=\phi(xt)$. 
    \end{enumerate}
   \indent Let $G$ be a countably infinite semigroup. A \textit{left Følner sequence} in $G$ is a sequence $(F_N)_{N=1}^\oo$ of finite non-empty subsets of $G$ satisfying
        \begin{equation} 
        \lim_{N\rightarrow \infty} \frac{|F_N \cap gF_N|}{|F_N|} = 1 
        \end{equation}
    for all $g \in G$. It is well-known that the existence of a left F\o lner sequence in a countably infinite semigroup $G$ implies that $G$ is left amenable (see \cite[Theorem 6.4]{F60}). Furthermore, if $G$ is cancellative,\footnote{A semigroup $G$ is (two-sided) \textit{cancellative} if, for any $a,b,c\in G$, $ab=ac\Longrightarrow b=c$ and $ba=ca\Longrightarrow b=c$.} then having a left F\o lner sequence is equivalent to $G$ being left amenable (see \cite[Theorem 2]{AW67}). From now on, we will tacitly assume that the semigroups we deal with are countable and cancellative.\\
    \indent Given a left Følner sequence $\mathbf{F} = (F_N)_{N=1}^\oo$ in an amenable semigroup $G$ and a set $S \sseq G$, define \begin{equation}
        \overline{d}_\mathbf{F}(S) = \limsup_{N\rightarrow \infty}\frac{|S \cap F_N|}{|F_N|}, \quad \quad \underline{d}_\mathbf{F}(S) = \liminf_{N\rightarrow \infty}\frac{|S \cap F_N|}{|F_N|}.
            \end{equation}
    If $\overline{d}_\mathbf{F}(S) = \underline{d}_\mathbf{F}(S)$, then their common value is denoted by $d_\mathbf{F}(S)$. \\
    \indent For the remainder of the paper we will, as a rule, omit the adjective \textit{left} when dealing with left amenable semigroups and left F\o lner sequences. \\
    \indent We will now embark on proving a version of Theorem \ref{Ruz} in the setup of amenable groups (Theorem \ref{general ruzsa} below). The further extension to amenable semigroups will be done in subsection \ref{Sec 4.3}.\\
    \indent Throughout the rest of this section, an important role is played by \textit{tempered} F\o lner sequences. 

    \begin{definition} \label{tempered} 
    For a group $G$, a sequence $(F_n)_{n=1}^\oo$ of finite subsets of $G$ is called \textit{tempered} if there is some constant $C>0$ such that  
            \begin{equation} \left| \bigcup_{k<n} F_k^{-1}F_n \right| \leq C|F_n| \end{equation}
        for all $n \in \N$ with $n>1$. (For $A,B\subseteq G$ we define $A^{-1}B\defeq\{a^{-1}b:a\in A,b\in B\}$).
        \end{definition}

In much of what we do in this section, the following pointwise ergodic theorem, which utilizes the notion of a tempered F\o lner sequence, is essential.

\begin{Theorem}[ (\cite{Shu88}, see Section 5.6 in \cite{TEM92})]  \label{Schulman}
    Let $G$ be a group with tempered Følner sequence $(F_N)_{N=1}^\infty$ and let $(X,\mathcal{B},\mu,(T_g)_{g \in G})$ be a measure-perserving system. Then for any $ f \in L^2$, the limit
        \[ \lim_{N \rightarrow \infty} \frac{1}{|F_N|} \sum_{g \in F_N} f(T_gx) \]
    exists a.e.
\end{Theorem}

In order to prove the generalization of Theorem \ref{Ruz} to countably infinite amenable groups we will need two results, Theorems \ref{ABC123} and \ref{general Furst}, which can be viewed, correspondingly, as amenable analogs of Theorems \ref{Berg} and \ref{Furst}.

    \begin{Theorem}\label{ABC123}
        Let $G$ be a countably infinite amenable group with a tempered Følner sequence $\mathbf{F} = (F_N)_{N=1}^\infty$. Let $(X, \mathcal{B}, \mu, (T_g)_{g \in G})$ be a measure-preserving system and let $B \in \mathcal{B}$ with $\mu(B) > 0$. Then there is a set $P \subseteq G$ such that
        \begin{equation}\label{DesiredResult}
            d_\bo{F}(P)\ge\mu(B)\quad\text{and}\quad \mu\left( B \cap T_{g_1}^{-1}B \cap \cdots \cap T_{g_r}^{-1}B \right)> 0\quad\text{for all}\quad 
            g_1,...,g_r\in P.
        \end{equation}
        
    \end{Theorem}


                \begin{Proof}
                 After deleting (if needed) a set of measure zero from each $T^{-1}_{g}B$, we can and will assume without loss of generality that for any $g_1,...,g_j\in G$, if $T^{-1}_{g_1} B\cap\cdots\cap T^{-1}_{g_j}B\neq\emptyset$, then $\mu(T^{-1}_{g_1} B\cap\cdots\cap T^{-1}_{g_j}B)>0$. For $N = 1,2,3,\ldots$ define $f_N:X \rightarrow \R$ by,
                    \[ f_N(x) = \frac{1}{|F_N|} \sum_{g \in F_N} 1_{B}({T_g}x) .\]
                For all $N\in\N$ we have $\int f_N = \mu(B)$. Let $f = \limsup_{N \rightarrow \infty} f_N$. Note that $0 \leq f_N \leq 1$ for all $N\in \N$, so $f$ is well-defined and takes values in $[0,1]$. Then by Fatou's lemma, we have
                    \[ \int f = \int \limsup_{N\rightarrow \infty} f_N \geq \limsup_{N\rightarrow \infty} \int f_N = \mu(B). \]
                Since $f$ is nonnegative and $\int f \geq \mu(B)$, the set $f^{-1}[\mu(B),1]$ must have positive measure. By Theorem \ref{Schulman}, $\lim_{N\to\oo} f_N$ exists a.e., so in particular there is some $x_0 \in f^{-1}[\mu(B), 1]$ such that $\lim_N f_N(x_0) = f(x_0)$. If we set $P' \defeq \{g \in G: T_{g} x_0 \in B\}$, then
                    \[ d_\mathbf{F}(P') = \lim_{N\to\oo} \frac{|F_N \cap P'|}{|F_N|} =\lim_{N\to \oo} \frac{1}{|F_N|} \sum_{g \in F_N} 1_{B}({T_g}x_0) = \lim_{N\to \oo} f_N(x_0) = f(x_0) \geq \mu(B). \]
                \indent Now choose $h \in G$ so that $T_hx_0 \in B$ and put $P = P'h^{-1}$. We claim that $P$ satisfies formula (\ref{DesiredResult}).
                Firstly, observe $d_\mathbf{F}(P) = d_\mathbf{F}(P') \geq \mu(B)$. Secondly, take $g_1h^{-1},\ldots,g_rh^{-1} \in P$, where $g_1,\ldots, g_r \in P'$. See that for $i = 1,\ldots,r$ we have $T_{g_i}x_0 \in B$, so $T_{g_ih^{-1}}T_hx_0 \in B$, therefore $T_hx_0 \in T_{g_ih^{-1}}^{-1}B$. Thus the intersection $B \cap T_{g_1h^{-1}}^{-1}B \cap \cdots \cap T_{g_rh^{-1}}^{-1}B $ is non-empty since it contains $T_hx_0$, so $\mu \left( B \cap T_{g_1h^{-1}}^{-1}B \cap \cdots \cap T_{g_rh^{-1}}^{-1}B\right)> 0$. 
                \end{Proof}

For the proof of Theorem \ref{general ruzsa} we will need the following variant of Furstenberg's correspondence principle 
in the amenable setup.

\begin{Theorem}[ (\site{Theorem 6.4.17}{Ber00})] \label{general Furst}
    Let $G$ be a countably infinite amenable group, let $\mathbf{F}$ be a Følner sequence in $G$, and let $S \sseq G$ with $\overline{d}_\mathbf{F}(S) > 0$. Then there is an invertible measure-preserving system $(X,\mathcal{B},\mu,(T_g)_{g \in G})$ and a set $B \in \mathcal{B}$ with $\mu(B) = \overline{d}_\mathbf{F}(S)$ such that
        \begin{equation}
            \overline{d}\left(S \cap g_1^{-1}S \cap \cdots \cap g_r^{-1}S \right)\geq\mu\left( B \cap T_{g_1}^{-1}B \cap \cdots \cap T_{g_r}^{-1}B \right) 
        \end{equation}
    for all $g_1,\ldots,g_r \in G$. 
\end{Theorem}

If $G$ is a semigroup and $A \sseq G$, then for all $g \in G$ define $g^{-1}A = \{h \in G: gh \in A\}$ and $Ag^{-1} = \{h \in G: hg \in A\}$.

\begin{definition} \label{delta sets in semigroups}
    For an amenable semigroup $G$ with Følner sequence $\mathbf{F}$ and $S \subseteq G$, let
        \begin{gather}
        \label{eq 4.10}
            \Delta_1(S)= \{g \in G : S \cap (g^{-1}S) \neq \emptyset\},  \\
        \label{eq 4.11}
            \Delta_2(\mathbf{F},S) = \{g \in G: \overline{d}_\mathbf{F}( S \cap (g^{-1}S)) >0\}.  
        \end{gather}
\end{definition}

\begin{remark} When $G$ is a group, the occurrences of $S \cap (g^{-1}S)$ in equations (\ref{eq 4.10}) and (\ref{eq 4.11}) can be replaced by $S \cap (gS)$ without changing the meaning of the Definition \ref{delta sets in semigroups}. In other words, when $G$ is a group, both $\Delta_1(S)$ and $\Delta_2(\mathbf{F},S)$ are symmetric, meaning that $g\in \Delta_1(S)\iff g^{-1}\in \Delta_1(S)$ and $g\in \Delta_2(\mathbf{F},S)\iff g^{-1}\in \Delta_2(\mathbf{F},S)$.\\
\indent When dealing with non-commutative groups there are two possible ways to define $\Delta_1(S)$ and $\Delta_2(\mathbf{F},S)$. For example, one could define $\Delta_1(S)$ as $   \{g \in G : S \cap (Sg^{-1}) \neq \emptyset\}$ which is equivalent to definition of $\Delta_1(S)$ in the commutative case.
The ``left" choice made in equations (\ref{eq 4.10}) and (\ref{eq 4.11}) is consistent with our definition of amenability via left invariant means and left F\o lner sequences. 
\end{remark}

We are now ready to state and prove a generalization of Theorem \ref{Ruz} to countably infinite amenable groups.



\begin{Theorem}[ (Theorem \ref{Ruz} for countably infinite amenable groups)] \label{general ruzsa}
    Let $G$ be a countably infinite amenable group with (left) Følner sequences $\mathbf{G},\mathbf{F}_1,\ldots,\mathbf{F}_k$, where $\mathbf{G}$ is (left) tempered. Let $S_1,\ldots,S_k \sseq G$ such that $\overline{d}_{\mathbf{F}_i}(S_i)> 0$ for each $i=1,\ldots,k$. Then there is $S \sseq G$ such that $d_\mathbf{G}(S) \geq \prod_{i=1}^k \overline{d}_{\mathbf{F}_i}(S_i)$ and $
                \Delta_1(S) \sseq D = \bigcap_{i=1}^k \Delta_2(\mathbf{F}_i,S_i)$. Furthermore, there are $m_1,\ldots,m_\ell \in G$, where
            \begin{equation}
                \ell \leq \prod_{i=1}^k 1/\overline{d}_{\mathbf{F}_i}(S_i),
            \end{equation}
        such that $\bigcup_{i=1}^\ell(m_iD) = \bigcup_{i=1}^\ell(Dm_i^{-1}) = G$. 
\end{Theorem}

    \begin{Proof}
   By Theorem \ref{general Furst}, for $i = 1,\ldots,k$ there is a measure-preserving system $(X_i, \mathcal{B}_i, \mu_i, (T_{i,g})_{g \in G})$ and there is $B_i \in \mathcal{B}_i$ with $\mu_i(B_i) = \overline{d}_{\mathbf{F}}(S_i)$ such that
                \[  \overline{d}_{\mathbf{F}_i} \left( g_1^{-1}S_i \cap \cdots \cap g_m^{-1}S_i \right) \geq \mu_i \left( T_{i,g_1}^{-1}B_i \cap \cdots \cap T_{i,g_m}^{-1}B_i \right) \]
            for all $g_1,\ldots,g_m \in G$. 
            For all $g \in G$ let $(X, \mathcal{B}, \mu, T_g) = \prod_{i=1}^k(X_i, \mathcal{B}_i, \mu_i, T_{i,g})$ and let $B = B_1 \times \cdots \times B_k$. By Lemma \ref{ABC123} there exists $S\subseteq G$ such that $d_{\mathbf{G}}(S)\geq \mu(B)$ and for all finite $g_1,\ldots,g_m \in S$ we have
            \[\mu\bigg(B \cap T_{g_1}^{-1}B \cap \cdots \cap T_{g_m}^{-1}B \bigg)>0.\]
            Note that
            \[\mu(B)=\prod_{i=1}^k\mu_i(B_i)=\prod_{i=1}^k \overline{d}_{\mathbf{F}_i}(S_i),\]
        so $d_{\mathbf{G}}(S) \geq \prod_{i=1}^k \overline{d}_{\mathbf{F}_i}(S_i)$. Recall $D$ is defined to be the set $\bigcap_{i=1}^k {\Delta_2(\mathbf{F}_i},S_i)$. In order to show that $\Delta_1 (S) \subseteq D$, take $s,t \in S$ and note that
            \begin{align*}
                0 < \mu(T_s^{-1} B \cap T_t^{-1} B) &= \mu( B \cap T_{ts^{-1}}^{-1} B) = \prod_{i=1}^k \mu_i( B_i \cap T_{i,ts^{-1}}^{-1} B_i) \leq \prod_{i=1}^k \overline{d}_{\mathbf{F}_i}(S_i \cap (ts^{-1})^{-1} S_i).
            \end{align*}
        It follows that $ts^{-1} \in \Delta_2(\mathbf{F},S_i)$ for all $i \in \{1,\ldots,k\}$, meaning that $ts^{-1} \in D$. Thus $\Delta_1(S) \subseteq D$, as desired. Lastly, there are $m_1,\ldots,m_\ell \in G$, where
            \[ \ell \leq \prod_{i=1}^k 1/\overline{d}_{\mathbf{F}_i}(S_i), \]
        such that $G = m_1\Delta_1(S) \cup \cdots \cup m_\ell \Delta_1(S)= \Delta_1(S) m_1^{-1} \cup \cdots \cup \Delta_1(S) m_\ell^{-1}$ (this will be proven below in Theorem \ref{4 D}). Since $\Delta_1(S)\subseteq D$ it follows that $G = m_1D \cup \cdots \cup m_\ell D= D m_1^{-1} \cup \cdots \cup D m_\ell^{-1}$. \\
        \hphantom{hi}
    \end{Proof}
        
\indent The set $S$ in Theorem 4.7 has the property that $\Delta_1(S)\subseteq\bigcap_{i=1}^k \Delta_2(\mathbf{F}_i,S_i)$ which means that for all $g$ such that $S\cap g^{-1}S\neq \emptyset$ it is true that $\overline{d}_{\mathbf{F}_i}(S_i \cap (g^{-1}S_i)) > 0$ for $1 \leq i \leq k$. In fact one can show $S$ satisfies the stronger property that for all $g_1,\dots,g_n\in G$ such that $S \cap (g_1^{-1}S) \cap \cdots \cap (g_n^{-1}S) \neq \emptyset$ it is true that $\overline{d}_{\mathbf{F}_i}(S_i \cap (g_1^{-1}S_i) \cap \cdots \cap (g_n^{-1}S_i)) > 0$ for $1 \leq i \leq k$. This result is reflected in the following theorem.

    \begin{Theorem} \label{4.9}
        Let $G$ be a countably infinite amenable group with (left) Følner sequences $\mathbf{G},\mathbf{F}_1,\ldots,\mathbf{F}_k$, where $\mathbf{G}$ is (left) tempered. Let $S_1,\ldots,S_k \sseq G$ such that $\overline{d}_{\mathbf{F}_i}(S_i)> 0$ for each $i=1,\ldots,k$. Then there is a set $S \sseq G$ such that $d_{\mathbf{G}}(S) \geq \prod_{i=1}^k \overline{d}_{\mathbf{F}_i}(S_i)$ and if $g_1,\ldots,g_n \in G$ are such that $S \cap (g_1^{-1}S) \cap \cdots \cap (g_n^{-1}S) \neq \emptyset$, then $\overline{d}_{\mathbf{F}_i}(S_i \cap (g_1^{-1}S_i) \cap \cdots \cap (g_n^{-1}S_i)) > 0$ for $1 \leq i \leq k$.\footnote{In the special case where $k=1$, $G=\Z$, and the F\o lner sequences $\mathbf{G}$ and $\mathbf{F}_1$ are the standard 
        F\o lner sequence $\{1,\ldots,N\}$, $N=1,2,3,...$, 
        Theorem \ref{4.9} reduces to a result of Ellis which was proven in \cite[Corollary 2.1.1]{Ber85} and \cite[Theorem 3.20]{Fur81}. It's worth mentioning that both proofs of Ellis' result rely on the pointwise ergodic theorem. Indeed, the proof of 
        \cite[Corollary 2.1.1]{Ber85} 
        relies on Theorem \ref{Berg} (see Remark \ref{2.2}) and the proof of \cite[Theorem 3.20]{Fur81}
        relies on \cite[Proposition 3.7]{Fur81}
        which establishes the existence of generic points for ergodic continuous maps of compact spaces with the help of the pointwise ergodic theorem.}
    \end{Theorem}


    Utilizing the fact that every Følner sequence has a tempered Følner subsequence (see \cite[Proposition 1.4]{Lin99}), it is straightforward to obtain the following variant of Theorem \ref{general ruzsa} with the slightly weaker assumption that the F\o lner sequence $\mathbf{G}$ may not be tempered. 

    \begin{Theorem}  \label{some theorem}
     Let $G$ be a countably infinite amenable group and let $\mathbf{G},\mathbf{F}_1,\ldots,\mathbf{F}_k$ be Følner sequences in $G$. If $S_1,\ldots,S_k \sseq G$ satisfy $\overline{d}_{\mathbf{F}_i}(S_i)> 0$ for each $i=1,\ldots,k$, then there is $S \sseq G$ such that $\ud_{\mathbf{G}}(S) \geq \prod_{i=1}^k \overline{d}_{\mathbf{F}_i}(S_i)$ and $
                \Delta_1(S) \sseq D = \bigcap_{i=1}^k \Delta_2(\mathbf{F}_i,S_i)$. Furthermore, there are $m_1,\ldots,m_\ell \in G$, where
            \begin{equation}
                \ell \leq \prod_{i=1}^k 1/\overline{d}_{\mathbf{F}_i}(S_i),
            \end{equation}
        such that $\bigcup_{i=1}^\ell(m_iD) = \bigcup_{i=1}^\ell(Dm_i^{-1}) = G$. 
\end{Theorem}


    To complete the logical circle, we still owe the reader the proof of Theorem \ref{4 D} which was used in the proof of Theorem \ref{general ruzsa}. Before proceeding, we need to introduce some preliminary definitions and prove some necessary lemmas. Let $G$ be a semigroup and let $S \subseteq G$. The set $S$ is \textit{left syndetic} if there are $c_1,\ldots,c_k \in G$ such that $c_1^{-1}S \cup \cdots \cup c_k^{-1}S = G$. Similarly $S$ is \textit{right syndetic} if there are $d_1,\dots,d_r\in G$ such that $Sd_1^{-1}\cup\dots\cup Sd_r^{-1}=G$. $S$ is \textit{right thick} if for every finite $F\subseteq G$, there exists $g\in G$ for which $Fg\subseteq S$.


The next three lemmas are standard. We supply the proofs for the convenience of the reader.

\begin{Lemma} \label{4 A}
     Let $G$ be an infinite group and let $T\subseteq G$ be right thick. Then there exists an injective sequence $(s_n)_{n=1}^\infty$ in  $T$ such that $s_i^{-1}s_j \in T$ whenever $i < j$. 
\end{Lemma}

\begin{Proof}
    If $T=G$, then the result is clear. So suppose that $T\neq G$ and fix $t\in G\setminus T$. We inductively construct an injective sequence $s_1,s_2,...$ in $G$ with our desired property. Pick arbitrary $s_1\in G$. Now let $n\in\N$ be arbitrary and suppose that we have chosen pairwise distinct $s_1,...,s_n\in G$ with $\{s_i^{-1}s_j: 1 \leq i < j\leq n \} \subseteq T$. Now $F=\set{s_1^{-1},...,s_n^{-1},ts_1^{-1},...,ts_n^{-1}}$ is finite, so take $s_{n+1}\in G$ such that $Fs_{n+1}\subseteq T$. Since $t\not\in T$, we have $s_{n+1}\neq s_j$ for any $1\le j\le n$. Furthermore one can now see that $\{s_i^{-1}s_j: 1 \leq i < j\leq n+1 \} \subseteq T$, so  by induction this gives us an infinite sequence $s_1,s_2,...$ with the desired properties.
\end{Proof} 


\begin{Lemma} \label{4 B}
    Let $G$ be an infinite group. If $S\subseteq G$ is not left syndetic, then $T=G\setminus S$ is right thick. 
\end{Lemma}

\begin{Proof}
     Let $F\subseteq G$ be finite and non-empty, say $F = \{g_1,\ldots,g_k\}$. Then $g_1^{-1}S \cup \cdots \cup g_k^{-1}S\neq G$ by non-syndeticity of $S$, so pick $a\in G\setminus (g_1^{-1}S \cup \cdots \cup g_k^{-1}S)$. It follows that for all $g \in F$ we have $ga \notin S$, hence $Fa\subseteq G \setminus S=T$, which implies that $T$ is right thick as desired. 
\end{Proof}


\begin{Lemma} \label{4 C}
    Let $G$ be a countably infinite amenable group and let $\mathbf{F}$ be a Følner sequence in $G$. Let $S \subseteq G$ with $\overline{d}_\mathbf{F}(S) > 0$. If $(s_n)_{n=1}^\infty$ is an injective sequence in $G$, then there are $i<j$ such that $s_i^{-1}s_j \in \Delta_1(S)$. 
\end{Lemma}

\begin{Proof}
    Let $(s_n)_{n=1}^\infty$ be an injective sequence in $G$. By the pigeonhole principle, there must be some $i< j $ for which $\overline{d}_\mathbf{F}((s_iS) \cap (s_jS)) > 0$. It follows that there are $u,v \in S$ such that $s_ju = s_iv$, thus $(s_i^{-1}s_j)u = v$, and therefore $u \in S \cap (s_i^{-1}s_j)^{-1} S$. It follows that $ S \cap (s_i^{-1}s_j)^{-1} S \neq \emptyset$, therefore $s_i^{-1}s_j \in \Delta_1(S)$. 
\end{Proof}

We can now prove Theorem \ref{4 D}, thereby completing the proof of Theorem \ref{general ruzsa}.

\begin{Theorem}
     \label{4 D}
    Let $G$ be a countably infinite amenable group, let $\mathbf{F}$ be a Følner sequence in $G$,  let $S \sseq G$ with $\overline{d}_\mathbf{F}(S) > 0$, and let $E = \Delta_1(S)$. Then there exist $b_1,...,b_\ell\in G$ with $\ell\le 1/\overline{d}_\mathbf{F}(S)$ for which $b_1E \cup \cdots \cup b_\ell E = E b_1^{-1} \cup \cdots \cup Eb_\ell^{-1}= G$ (that is, $E$ is both left and right syndetic with the same ``syndeticity constant" $\ell$). 
\end{Theorem}

\begin{Proof}
 By Lemma \ref{4 C}, for every injective sequence $(s_n)_{n=1}^\infty$ in $G$ there are $i<j$ such that $s_i^{-1}s_j \in E$. Then by Lemma \ref{4 A} it follows that $G \setminus E$ cannot be right thick, and so $E$ is left syndetic by Lemma \ref{4 B},  which implies that there are $c_1,...,c_k\in G$ such that $c_1E\cup\cdots\cup c_kE=G$. \\
\indent For any additional $g\in G$ with $g \notin \{c_1,\ldots,c_k\}$ we have $g \in G = c_1E\cup\cdots\cup c_kE$, so $g\in c_jE$ for some $j\in\set{1,...,k}$. Hence $c_j^{-1}g \in E$, which means that $S \cap (c_j^{-1}g)^{-1}S \neq \emptyset$, so there are $u,v \in S$ such that $c_j^{-1}gu = v$. Thus  $c_j\{u,v\} \cap g\{u,v\} \neq \emptyset$. It follows that there is a maximal subset $\set{b_1,...,b_\ell}$ in $G$ such that, for all finite $F \subseteq S$, we have $b_iF\cap b_jF=\emptyset$ when $i \neq j$, and if $g \notin \{b_1,\ldots,b_\ell\}$, then there is some finite $F \sseq S$ such that $b_iF \cap g F\neq \emptyset$ for some $1 \leq i \leq \ell$.\\
\indent Write $\mathbf{F}=(F_n)_{n=1}^\oo$ and for all $n\in\N$ put $S(n) = |A_n|$, where $A_n\defeq S\cap F_n$. By assumption, for all $n\in\N$, the sets $b_1A_n,...,b_\ell A_n$ are disjoint. Hence,
\equat{
    \ell S(n)=\abs{b_1A_n\cup\cdots\cup b_\ell A_n}\le\abs{b_1F_n\cup\cdots\cup b_\ell F_n}.
}
Thus,
\begin{equation}
\label{eq 100000}
    \ell\le \frac{|b_1F_n\cup\cdots\cup b_\ell F_n|}{|F_n|}\cdot\frac{|F_n|}{S(n)}.
\end{equation}

By taking $n\to\oo$ in formula (\ref{eq 100000}) and using the fact that $\bo{F}$ is a F\o lner sequence, it follows that $\ell\le 1/\overline{d}_\mathbf{F}(S)$. \\
\indent We now show that $ b_1E\cup\cdots\cup  b_\ell E=G$. Indeed, let $g\in G$ be arbitrary. By maximality of the set $\{b_1,\ldots,b_\ell\}$, pick $i\in\set{1,...,\ell}$ and finite $F\subseteq S$ with $b_iF\cap gF\neq\emptyset$. Then $b_it_1=gt_2$ for some $t_1,t_2\in F$. Since $(t_1t_2^{-1})t_2= t_1$ it follows that $S \cap (t_1t_2^{-1})S \neq \emptyset$, so $t_1t_2^{-1} \in D$, and since $b_it_1 = gt_2$ it follows that $b_i^{-1}g = t_1t_2^{-1} \in E$, thus $g \in b_iE$. This shows that $G=b_1E\cup\cdots\cup b_\ell E$.\\
\indent Lastly, note that if $d \in E$, then $d^{-1} \in E$. For all $g \in G$ we have $g^{-1} \in b_1E \cup \cdots \cup b_\ell E$, hence $b_i^{-1}g^{-1} = d$ for some $d \in E$ and so $gb_i = d^{-1} \in E$, thus $g \in Eb_i^{-1}$. This shows that $G = Eb_1^{-1} \cup \cdots \cup Eb_\ell^{-1}$. 
\end{Proof}

    \subsection{Extending Theorem \ref{general ruzsa} to Cancellative Amenable Semigroups } \label{Sec 4.3}

  
        \indent In this subsection we will obtain a variant of Theorem \ref{general ruzsa} for cancellative amenable semigroups, namely Theorem \ref{general ruzsa for semigroups}. We will do so by utilizing the classical fact that cancellative amenable semigroups can be embedded into groups. Let $G$ be a cancellative amenable semigroup. It is known that there is a group $\widetilde{G}$ and an embedding $\varphi\col G\to\widetilde{G}$ such that $\widetilde{G}=\set{\varphi(s)\varphi(t)^{-1}:s,t\in G}$ (see \cite[Theorem 1.23]{CP61} and \cite[Proposition 1.23]{PAT88}). Moreover, the group $\widetilde{G}$ is 
        unique (see \cite[Theorem 1.25]{CP61} for the precise formulation). With this in mind, we henceforth, when convenient, identify every element $g\in G$ with $\varphi(g)\in\widetilde{G}$ and call $\widetilde{G}$ the \textit{group of quotients} of $G$.\\
    \indent The following theorem allows us to relate F\o lner sequences in $G$ to F\o lner sequences in $\widetilde{G}$.


\begin{Theorem}[ (\site{Theorem 2.12}{BDM20})]
\label{Group Embedd}
    Let $G$ be a countably infinite amenable cancellative semigroup. Then any F\o lner sequence $(F_N)_{N=1}^\infty$ in $G$ is a F\o lner sequence in its group of quotients $\widetilde{G}$.
\end{Theorem}

\begin{remark}
\label{4.15}
    Let $G$ be a countably infinite amenable cancellative semigroup and let $\mathbf{F}=(F_N)_{N=1}^\oo$ be a  F\o lner sequence in $G$, and hence in $\widetilde{G}$.
    Let $S'\subseteq\widetilde{G}$ and let $S\defeq S'\cap G$. We have
    \[\ud_{\mathbf{F}}(S)=\limsup_{N\to\oo}\frac{|S\cap F_N|}{|F_N|}=\limsup_{N\to\oo}\frac{|S'\cap F_N|}{|F_N|}=\ud_{\mathbf{F}}(S').\]
    Moreover, if one of $d_{\mathbf{F}}(S)$ or $d_{\mathbf{F}}(S')$ exists, then both exist and $d_{\mathbf{F}}(S)=d_{\mathbf{F}}(S')$.
\end{remark}

We will need one additional fact about $\widetilde{G}$ which is encompassed in the following lemma.
    
    \begin{Lemma}
        \label{Coroll Thick}
        If $G$ is a countably infinite (left) amenable cancellative semigroup and $\widetilde{G}$ is its group of quotients, then for any finite $F\subseteq \widetilde{G}$, there exists $g\in G$ such that $Fg\subseteq G$. In particular, $G$ is right thick in $\widetilde{G}$.
    \end{Lemma}

        \begin{Proof}
            Let $(F_N)_{N=1}^\oo$ be a (left) F\o lner sequence in $G$, and hence a F\o lner sequence in $\widetilde{G}$ by Theorem \ref{Group Embedd}. Let $F=\{g_1,\dots, g_k\}$ be any finite subset of $\widetilde{G}$. Since $(F_N)_{N=1}^\oo$ is a F\o lner sequence in $\widetilde{G}$,
            \[\lim_{N\rightarrow \infty} \frac{\left|F_N \cap g_1^{-1}F_N\cap g_2^{-1}F_N\cap\dots \cap g_n^{-1}F_N\right|}{|F_N|}=1.\]
            This means there exists $r\in\N$ such that $F_r\cap g_1^{-1}F_r\cap g_2^{-1}F_r\cap\dots \cap g_k^{-1}F_r\neq \emptyset$. Now, for any element $g\in F_r \cap g_1^{-1}F_r\cap g_2^{-1}F_r\cap\dots \cap g_k^{-1}F_r$, we have $Fg\subseteq F_r\subseteq G$. Since $F$ was an arbitrary finite subset of $\widetilde{G}$, it follows that $\bigcup_{N=1}^\oo F_N$ is right thick in $\widetilde{G}$. Since $(F_N)_{N=1}^\oo$ is a F\o lner sequence in $G$, it follows that $\bigcup_{N=1}^\oo F_N\subseteq G$, and hence $G$ is right thick in $\widetilde{G}$.
        \end{Proof}

        \begin{remark}
            If $S\subseteq G$ is right thick in $G$, then $S$ is right thick in $\widetilde{G}$. Indeed, for any finite subset $F$ of $\widetilde{G}$ there exists a $g_1\in \widetilde{G}$ such that $Fg_1\subseteq G$ by Lemma \ref{Coroll Thick}. Since $Fg_1$ is a finite subset of $G$ there exists $g_2\in G$ such that $Fg_1g_2\subseteq S$. 
        \end{remark}


To obtain a variant of Theorem \ref{general ruzsa} for semigroups, one needs to slightly modify the definition of tempered F\o lner sequences, but first we will need to introduce some notation.
Let $G$ be a semigroup. For $B \sseq G$, define
            \begin{equation}\label{A-B_def}
                A^{-1}B = \bigcup_{a \in A}a^{-1}B,
            \end{equation} 
 where $a^{-1}B=\{h \in G: ah \in B\}$. Note that when $G$ is a group, $A^{-1}B=\{a^{-1}b:a\in A,b\in B\}$.

 \begin{definition} \label{tempered in semigroups}
     A sequence $(F_n)_{n=1}^\oo$ of finite subsets of a semigroup $G$ is \textit{(left) tempered} if for some $C>0$
        \begin{equation}
        \label{4.11}
            \left| \bigcup_{k < n} F_k^{-1}(F_ng)\right| \leq C|F_n|. 
        \end{equation} 
        for all $n \in \N$ with $n>1$, and for all $g\in G$.
 \end{definition}
 
 \begin{remark}
    When $G$ is a group, it can be seen that for all $n\in\N$ with $n>1$, and for all $g \in G$,
        \[\left| \bigcup_{k < n} F_k^{-1}(F_ng)\right|=\left| \bigcup_{k < n} F_k^{-1}F_n\right|.\]
        In this case Definition \ref{tempered in semigroups} agrees with Definition \ref{tempered}. 
        Note that for semigroups one may have that $(F_k^{-1}F_n)g\neq F_k^{-1}(F_ng)$.
        The choice of parentheses which is needed when proving Theorem \ref{tempered iff condition} is $F_k^{-1}(F_ng)$, for this reason we will write $F_k^{-1}F_ng\defeq F_k^{-1}(F_ng)$.
 \end{remark}



\begin{Theorem} \label{tempered iff condition}
    Let $G$ be a countably infinite amenable cancellative semigroup and $\widetilde{G}$ be its group of quotients. Let $(F_n)_{n=1}^\oo$ be a F\o lner sequence in $G$. Then the sequence $(F_n)_{n=1}^\oo$ is a tempered F\o lner sequence in $G$ iff it is a tempered F\o lner sequence in $\widetilde{G}$.
\end{Theorem}

\begin{Proof} The fact that $(F_n)_{n=1}^\infty$ is a Følner sequence in $G$ iff it is a Følner sequence in $\widetilde{G}$ follows from Theorem \ref{Group Embedd}. Let $\varphi\col G\to\widetilde{G}$ be an embedding of $G$ into its group of quotients $\widetilde{G}$. We wish to show that the sequence $(F_n)_{n=1}^\infty$ is a tempered F\o lner sequence in $G$ if and only if $(\varphi(F_n))_{n=1}^\infty$ is a tempered F\o lner sequence in $\widetilde{G}$.
Note that for all $a\in G$ and for $n\in\N$ with $n>1$,
    \[\varphi\left(\bigcup_{k=1}^{n-1}F_k^{-1} F_n a\right)\subseteq \bigcup_{k=1}^{n-1}\varphi(F_k)^{-1} \varphi(F_n)\varphi(a).\]
    Now suppose that $(\varphi(F_n))_{n=1}^\infty$ is a tempered F\o lner sequence in $\widetilde{G}$. Then there exists $C>0$ such that for all $a\in G$ and for $n\in\N$ with $n>1$,
    \begin{align*}
        \left|\bigcup_{k=1}^{n-1}F_k^{-1} F_na\right|=\left|\varphi\left(\bigcup_{k=1}^{n-1}F_k^{-1} F_na\right)\right|&\leq\left| \bigcup_{k=1}^{n-1}\varphi(F_k)^{-1}\varphi(F_n) \varphi(a)\right|\\
        &=\left|\bigcup_{k=1}^{n-1}\varphi(F_k)^{-1} \varphi(F_n)\right|\leq C|\varphi(F_n)|=C|F_n|.
    \end{align*}
    This shows that $(F_n)_{n=1}^\oo$ is a tempered F\o lner sequence in $G$. Now assume that $(\varphi(F_n))_{n=1}^\infty$ is not a tempered F\o lner sequence in $\widetilde{G}$. This means that for all $C>0$ there exists $n\in \N$ with $n>1$ such that
    \[\left|\bigcup_{k=1}^{n-1}\varphi(F_k)^{-1} \varphi(F_n)\right|> C|\varphi(F_n)|.\]
    Now note that, for $n\in \N$ with $n>1$, since $\bigcup_{k=1}^{n-1}\varphi(F_k)^{-1} \varphi(F_n)$ is a finite subset of $\widetilde{G}$, by Lemma \ref{Coroll Thick} there exists an $a_n\in G$ such that $\bigcup_{k=1}^{n-1}\varphi(F_k)^{-1} \varphi(F_n)\varphi(a_n)\subseteq \varphi(G)$. This means 
    \[\bigcup_{k=1}^{n-1}\varphi(F_k)^{-1} \varphi(F_n)\varphi(a_n)=\varphi\left(\bigcup_{k=1}^{n-1}F_k^{-1} F_n a_n\right).\]
    So for all $C>0$ there exists $n\in \N$ with $n>1$ and $a_n\in G$ such that
    \begin{align*}
        \left|\bigcup_{k=1}^{n-1}F_k^{-1} F_n a_n\right|=\left|\varphi\left(\bigcup_{k=1}^{n-1}F_k^{-1} F_na_n\right)\right|&= \left|\bigcup_{k=1}^{n-1}\varphi(F_k)^{-1} \varphi(F_n)\varphi(a_n)\right|\\
    &=\left|\bigcup_{k=1}^{n-1}\varphi(F_k)^{-1} \varphi(F_n)\right|> C|\varphi(F_n)|=C|F_n|
    \end{align*}
    This shows that $(F_n)_{n=1}^\oo$ is not a tempered F\o lner sequence in $G$.
\end{Proof}

    Before proceeding to the proof Theorem \ref{general ruzsa for semigroups}, we need one final lemma. 

\begin{Lemma} \label{syndeticity constant for semigroups}
    Let $G$ be a countably infinite amenable cancellative semigroup, let $\mathbf{F}$ be a Følner sequence in $G$,  let $S \sseq G$ with $\overline{d}_\mathbf{F}(S) > 0$, and let $E = \Delta_1(S)$. Then there exist $b_1,\ldots,b_\ell \in G$ with $\ell \leq 1/\overline{d}_\mathbf{F}(S)$ for which $E b_1^{-1} \cup \cdots \cup Eb_\ell^{-1}= G$. 
\end{Lemma}

    \begin{Proof} Let $\widetilde{E}=\{g\in\widetilde{G}: S\cap g^{-1}S\neq\emptyset\}$. Now, since $\bo{F}$ is a F\o lner sequence in $\widetilde{G}$ by Theorem \ref{Group Embedd}, it follows from Theorem \ref{4 D} that there exists $k\le 1/\ud_\bo{F}(S)$ (here, we regard $S$ as a subset of $\widetilde{G}$) and $a_1,...,a_k\in\widetilde{G}$ such that $\widetilde{E}a_1^{-1}\cup\cdots\cup \widetilde{E}a_k^{-1}=\widetilde{G}$.
    A look at the proof of Theorem \ref{4 D} reveals that the set $\{a_1,\ldots,a_k\}$ has an additional property: $a_iF \cap a_jF =\emptyset$ for all finite $F \sseq S$ and all $i \neq j$. In fact $\{a_1,\ldots,a_k\}$ is maximal with respect to this property in the sense that for any set $A \sseq \widetilde{G}$ with more than $k$ elements there are $a,b \in A$ such that $aF \cap bF \neq \emptyset$ for some finite $F \sseq S$. Therefore one may choose a subset $\{b_1,\ldots,b_\ell\} $ of $G$, where $\ell \leq k$, such that if $g\in \widetilde{G}\setminus\set{b_1,...,b_\ell}$, then there is some finite set $F \sseq S$ and $1 \leq i \leq \ell$ such that $gF \cap b_iF \neq \emptyset$. Note that $\ell \leq k\leq 1/\overline{d}_\mathbf{F}(S)$. \\
    \indent We now show that $Eb_1^{-1} \cup \cdots \cup Eb_\ell^{-1} = G$, as desired. Let $g \in G$. Then $g$ has some inverse $g^{-1}$ in the group $\widetilde{G}$. By the property of the set $\{b_1,\ldots,b_\ell\}$ there is a finite set $F \sseq S$ such that $g^{-1}F \cap b_iF \neq \emptyset$ for some $1 \leq i \leq \ell$,  so there exist $u,v\in F$ such that $gb_iv=u$. We have $v \in S$ and $gb_iv = u \in S$, so $v \in S \cap (gb_i)^{-1}S$, thus $gb_i \in E$, meaning $g \in Eb_i^{-1}$.  \end{Proof}

With the preliminary results at hand, we are now ready to prove a variant of Theorem \ref{general ruzsa} for countably infinite amenable cancellative semigroups.

\begin{Theorem}
    \label{general ruzsa for semigroups}
    Let $G$ be a countably infinite amenable cancellative semigroup and let $\mathbf{G},\mathbf{F}_1,\ldots,\mathbf{F}_k$ be Følner sequences in $G$ such that $\mathbf{G}$ is tempered. If $S_1,\ldots,S_k \sseq G$ satisfy $\overline{d}_{\mathbf{F}_i}(S_i)> 0$ for each $i=1,\ldots,k$, then there is $S \sseq G$ such that $d_\mathbf{G}(S) \geq \prod_{i=1}^k \overline{d}_{\mathbf{F}_i}(S_i)$ and $
                \Delta_1(S) \sseq D := \bigcap_{i=1}^k \Delta_2(\mathbf{F}_i,S_i)$. Furthermore, there are $m_1,\ldots,m_\ell \in G$, where
            \begin{equation}
                \ell \leq \prod_{i=1}^k 1/\overline{d}_{\mathbf{F}_i}(S_i),
            \end{equation}
        such that $\bigcup_{i=1}^\ell(Dm_i^{-1}) = G$. If in addition $G$ is a group, then $\bigcup_{i=1}^\ell(m_iD) = G$ as well. 
\end{Theorem}

\begin{Proof}
     Let $\widetilde{G}$ be the group of quotients of $G$. Since in this proof we will be working simultaneously with subsets of $G$ and $\widetilde{G}$, it will be convenient to slightly modify the notation introduced in Definition \ref{delta sets in semigroups}. Namely for any $A\subseteq G$, $B\subseteq\widetilde{G}$, and Følner sequence $\mathbf{F}$ in $G$ (considered, when convenient, as a F\o lner sequence in $\widetilde{G}$) we will write
    \[\Delta_1(A,G) = \{g \in G: A \cap g^{-1}A \neq \emptyset\} \quad \text{and} \quad \Delta_2(\mathbf{F},A,G) = \{ g \in G: \overline{d}_\mathbf{F}(A \cap g^{-1}A) > 0\}, \]
        \[ \Delta_1(B,\widetilde{G}) = \{g \in \widetilde{G}: B \cap g^{-1}B \neq \emptyset\} \quad \text{and} \quad \Delta_2(\mathbf{F},B,\widetilde{G}) = \{g \in \widetilde{G}: \overline{d}_\mathbf{F}(B \cap g^{-1}B) > 0\}. \]
    Theorem \ref{Group Embedd} implies that $\mathbf{G},\mathbf{F}_1,\ldots,\mathbf{F}_k$ are F\o lner sequences in $\widetilde{G}$. Additionally, Theorem \ref{tempered iff condition} shows $\mathbf{G}$ is tempered in $\widetilde{G}$. So by Theorem $\ref{general ruzsa}$, there exists $S'\subseteq \widetilde{G}$ such that $d_\mathbf{G}(S') \geq \prod_{i=1}^k \overline{d}_{\mathbf{F}_i}(S_i)$ and $\Delta_1(S',\widetilde{G})\subseteq \bigcap_{i=1}^k \Delta_2(\mathbf{F}_i,S_i,\widetilde{G})$. Let $S\defeq S'\cap G$. Note that by Remark \ref{4.15} 
    \[d_{\mathbf{G}}(S)=d_\mathbf{G}(S') \geq \prod_{i=1}^k  \overline{d}_{\mathbf{F}_i}(S_i). \]
    Also observe that,
    \[\Delta_1(S,G)\subseteq \Delta_1(S',G)\subseteq \Delta_1(S',\widetilde{G})\subseteq \bigcap_{i=1}^k \Delta_2(\mathbf{F}_i,S_i,\widetilde{G}) \]
    which implies
    \[ \Delta_1(S,G)\subseteq \bigcap_{i=1}^k \Delta_2(\mathbf{F}_i,S_i,\widetilde{G})\cap G=D.\]
   Lastly, since $\Delta_1(S,G) \sseq D$, by Lemma \ref{syndeticity constant for semigroups} there are $m_1,\ldots,m_\ell \in G$ with 
    \[ \ell \leq \prod_{i=1}^k 1/\overline{d}_{\mathbf{F}_i}(S_i)\]
    such that $(Dm_1^{-1}) \cup \cdots \cup (Dm_k^{-1}) = G$. If $G$ is a group, then the fact that $\bigcup_{i=1}^\ell (m_i D) = G$ follows from Theorem \ref{general ruzsa}. 
\end{Proof}

\section{Examples of tempered Følner sequences}
In this section we collect some examples of tempered F\o lner sequences, which, as we saw in subsection \ref{Sec 4.2}, play an important role in extending Theorem \ref{Ruz} to the framework of amenable groups. In subsection \ref{Sec 5.1} we present some useful facts and examples relating to tempered and non-tempered Følner sequences in $(\N,+)$. In subsection \ref{Sec 5.2} we provide examples of tempered F\o lner sequences in $(\N,\cdot)$. In subsection \ref{Sec 5.3} we provide an example of a tempered F\o lner sequence in the Heisenberg group. Finally, in subsection \ref{Sec 5.4} we observe that the natural F\o lner sequences in locally finite groups are actually tempered. 
\subsection{Tempered Sequences in $(\N,+)$}\label{Sec 5.1}

\indent It is well-known that the sequence $(F_N)_{N=1}^\infty$, where $F_N = \{1,\ldots,N\}$,  is a tempered Følner sequence in $(\N,+)$. The following simple observation is a generalization of this fact:

\begin{Theorem}
\label{temp in N}
    For $a,b \in \N$, $a \leq b$, let $[a,b]$ denote the set $\{a,a+1,\dots,b-1,b\}$. 
    If $(a_n)_{n=1}^\infty$ and $(b_n)_{n=1}^\infty$ are sequences in $\N$ such that $b_n-a_n\to\infty$ and $[a_j,b_j]\subseteq [a_{j+1},b_{j+1}]$ for all $j\in\N$, then $F_n:=[a_n,b_n]$, $n=1,2,3,...$ is a tempered F\o lner sequence in $\N$.
\end{Theorem}

\begin{Proof}
  In view of Theorem \ref{tempered iff condition}, it suffices to show that $(F_n)$ is a tempered Følner sequence in $(\Z,+)$.
 For this reason, any differences $F_i - F_j$ will be computed in $(\Z,+)$.\footnote{Given $A,B\subseteq\Z$, we define $A-B=\set{a-b:a\in A,b\in B}.$ Note that this is the additive analog of formula (\ref{A-B_def}).} It is clear that $(F_n)$ is a F\o lner sequence in $(\Z,+)$. 
    Since $F_1\subseteq F_2\subseteq F_3\subseteq\cdots$, it follows that for any $n\in\N$ with $n>1$, we have
     \begin{align*}
        \left| \bigcup_{1 \leq k < n} (F_n-F_k)\right|&
        = |F_n - F_{n-1}|= \Big|[a_n-b_{n-1},b_n-a_{n-1}]\Big|\\
        &=b_{n}-a_{n-1}+b_{n-1}-a_n+1 \leq 2b_n-2a_n+1< 2(b_n-a_n+1) = 2|F_n|.
    \end{align*}    
    This shows that $(F_n)$ is tempered.
\end{Proof}




In \cite{ADJ75} it was shown that the pointwise ergodic theorem fails along the sequence $([n,n+[\sqrt{n}]])_{n=1}^\infty$, so it follows from Theorem \ref{Schulman} that the F\o lner sequence $([n,n+[\sqrt{n}]])_{n=1}^\infty$ is not tempered. This 
result was later strengthened in \cite[Theorem 1]{JS77} and \cite[Corollary 3.9]{DJR79}, where it was shown that for any nondecreasing sequence $(b_n)_{n=1}
^\oo$ in $\N$ such that $\lim_{n\to\oo}n/b_n=0$, the pointwise ergodic theorem fails along the F\o lner sequence $([n,n+b_n])_{n=1}^\oo$, and hence $([n,n+b_n])_{n=1}^\oo$ is not tempered. Another example of a F\o lner sequence that is not tempered is the 
sequence $([n^2,n^2+n])_{n=1}^\infty$ (see \cite[Corollary 6]{BJR90} and also \cite[Example 2.7]{RW92}). Since every F\o lner sequence has a tempered subsequence, it is natural to ask for an example of a tempered subsequence of $([n^2,n^2+n])_{n=1}
^\infty$ that is tempered. One such example was provided in \cite[page 53]{BJR90}, where it was shown that $([2^{2^n},2^{2^n}+2^{2^{n-1}}])_{n=1}^\infty$ is tempered subsequence of $([n^2,n^2+n])_{n=1}^\infty$.



\subsection{Tempered Sequences in $(\N,\cdot)$}\label{Sec 5.2}

In what follows, $(\N,\cdot)$ denotes the semigroup of the natural numbers under multiplication, and $\Q_{>0}^\times$ denotes the group of positive rational numbers under multiplication. Let $p_1 < p_2 < p_3 < \cdots$ denote the sequence of primes. 


\begin{Theorem} \label{Tempered non-isotropic}
    For all $n \in \N$, let
    \[F_n=\{ p_1^{c_1} \dots p_n^{c_n}:0\leq c_i\leq (i+1)^{2n}\}.\]
    Then $(F_n)_{n=1}^\oo$ is a tempered Følner sequence in $(\N,\cdot)$.
\end{Theorem}

\begin{Proof}
    We first show that $(F_n)$ is a F\o lner sequence. To this end, let $g\in \N$  have prime factorization $g=p_1^{r_1}\dots p_k^{r_k}$, where $r_j\geq 0$ for $1\leq j\leq k$, and observe that for all $n>k$,
    \[F_n\cap gF_n=\left\{p_1^{\ell_1}\dots p_{n}^{\ell_n}: r_i\leq\ell_i\leq(i+1)^{2n}\text{ for } 1\leq i\leq k,\text{ and } 0\leq\ell_i\leq(i+1)^{2n}\text{ for }i>k\right\}.\]
    This means that, for sufficiently large $n$,
    \[|F_n\cap gF_n|=\prod_{i=1}^k((i+1)^{2n}-r_i+1)\prod_{i=k+1}^n((i+1)^{2n}+1).\]
     Since $|F_n|=\prod_{i=1}^n ((i+1)^{2n}+1)$ for all $n\in\N$, we have
    \[\lim_{n\to\oo}\frac{|F_n\cap gF_n|}{|F_n|}
    =\lim_{n\to\oo}\prod_{i=1}^k\frac{(i+1)^{2n}+1-r_i}{(i+1)^{2n}+1}
    =1,\]
    and so $(F_n)_{n=1}^\oo$ is a F\o lner sequence in $(\N,\cdot)$. It follows from Theorem \ref{Group Embedd} that $(F_n)_{n=1}^\infty$ is also a Følner sequence in $\Q_{>0}^\times$.  Thus, by Theorem \ref{tempered iff condition}, it suffices to show that $(F_n)_{n=1}^\infty$ is tempered in $\Q_{>0}^\times$. For this reason, the remainder of our calculations will be carried out in the group $\Q_{>0}^\times$. 
    Let $n\in\N$ with $n>1$ and observe that since $F_1\subseteq F_2\subseteq F_3\subseteq\cdots$, we have
    \equat{
        \bigcup_{k=1}^{n-1}F_k^{-1}&F_n=F_{n-1}^{-1}F_n\\
        &=
        \set{
            p_1^{\ell_1}\dots p_n^{\ell_n}:\hspace{-0.1cm}
            \begin{array}{c}
                 -(i+1)^{2(n-1)}\leq \ell_i\leq (i+1)^{2n}\hspace{0.15cm}\text{for}\hspace{0.15cm}
                 1\leq i\leq n-1,\,  0\leq \ell_n\leq (n+1)^{2n}
            \end{array}\hspace{-0.1cm}
        }.
    }
    Therefore,
    \[\frac{|F_{n-1}^{-1}F_n|}{|F_n|}=\frac{((n+1)^{2n}+1)\prod_{i=1}^{n-1}((i+1)^{2(n-1)}+(i+1)^{2n}+1)}{\prod_{i=1}^n((i+1)^{2n}+1)}=\prod_{i=1}^{n-1}\left(1+\frac{(i+1)^{2(n-1)}}{(i+1)^{2n}+1}\right)\]
    \[\leq \prod_{i=1}^{n-1}\left(1+\frac{(i+1)^{2(n-1)}}{(i+1)^{2n}}\right)=\prod_{i=1}^{n-1}\left(1+\frac{1}{(i+1)^2}\right)\leq \prod_{i=1}^\oo\left(1+\frac{1}{(i+1)^2}\right).\]
     Let $C=\prod_{i=1}^\oo\left(1+\frac{1}{(i+1)^2}\right)$. Notice $C<\oo$. We have that for all $n>1$,
    \[\left|\bigcup_{k=1}^{n-1}F_k^{-1}F_n\right|=|F_{n-1}^{-1}F_n|\leq C|F_n|,\]
    and hence $(F_n)_{n=1}^\oo$ is a tempered F\o lner sequence in $(\N,\cdot)$.
\end{Proof}

\begin{remark}
    Note that for any $\epsilon>0$ if we let $F_n=\{ p_1^{c_1} \dots p_n^{c_n}:0\leq c_i\leq (i+1)^{(1+\epsilon)n}\}$ for all $n \in \N$, then an argument similar to the proof of Theorem \ref{Tempered non-isotropic} shows that $(F_n)_{n=1}^\oo$ is tempered. On the other hand, if $G_n=\{ p_1^{c_1} \dots p_n^{c_n}:0\leq c_i\leq (i+1)^{n}\}$ for all $n \in \N$, one can show that $(G_n)_{n=1}^\infty$ is not tempered. 
\end{remark}

    \begin{Theorem} \label{thm 5.3}
        Let $f: \N \rightarrow \N$ be a nondecreasing function such that $f(n) \rightarrow \infty$. For $n \in \N$ let $ F_n = \{p_1^{c_1} \cdots p_n^{c_n} : 0 \leq c_i \leq f(n) \text{ for } 1 \leq i \leq n \} $.
        Then $(F_n)_{n=1}^\infty$ is a Følner sequence in $(\N,\cdot)$ and it is tempered iff the sequence 
            \begin{equation} \label{seq iff BB}  \frac{nf(n)}{f(n+1)}, \quad n \in \N \end{equation}
        is bounded. 
    \end{Theorem}

    \begin{Proof} It is straightforward to verify that $(F_n)$ is a Følner sequence in $(\N,\cdot)$, so we only show that $(F_n)$ is tempered iff the sequence in formula (\ref{seq iff BB}) is bounded.
        By Theorem \ref{tempered iff condition}, it suffices to show that $(F_n)$ is a tempered Følner sequence in $\Q_{>0}^\times$ iff the sequence in formula (\ref{seq iff BB}) is bounded, so we will carry out all of our calculations in $\Q^\times_{>0}$ for the remainder of the proof.  Since $f$ is non-decreasing, for all $n\in\N$ with $n>1$, we have
    \equat{
        \frac{\abs{\bigcup_{k=1}^{n-1}F_{k}^{-1}F_n}}{\abs{F_n}}&=
        \frac{\abs{F_{n-1}^{-1}F_n}}{\abs{F_n}}=
        \frac{
            \abs{
            \set{p_1^{\ell_1}\cdots p_n^{\ell_n}:-f(n-1)\le \ell_1,\ell_2,...,\ell_{n-1}\le f(n), 0\le \ell_n\le f(n)}
            }
        }
        {\abs{F_n}}\\
        &
            =\frac{(f(n) + 1) (f(n) + f(n-1) + 1)^{n-1}}{(f(n) + 1)^n}=\leri{1 + \frac{f(n-1)}{f(n) + 1}}^{n-1}.
    }
    Hence, by the definition of temperedness, $(F_n)_{n=1}^\oo$ is tempered iff the sequence 
    \begin{equation} \label{seq iff AA} 
    a_n:=\leri{ 1 + \frac{f(n-1)}{f(n) + 1}}^{n-1}, \quad n \in \N  \end{equation}
    is bounded. \\
    \indent Now note that the sequence in formula (\ref{seq iff BB}) is bounded iff the sequence $b_n:=\frac{(n-1)f(n-1)}{f(n) + 1}$, $n \in \N$  is bounded. Thus, it suffices to show that $(a_n)_{n=1}^\infty$ is bounded iff the  $(b_n)_{n=1}^\infty$
    is bounded. \\
    \indent $(\Rightarrow)$ Suppose that $(a_n)_{n=1}^\infty$ is bounded. Since for every $n>1$, we have 
    \equat{
        b_n=\frac{(n-1)f(n-1)}{f(n)+1}\le
        \sum_{k=0}^{n-1}\binom{n-1}{k}\leri{\frac{f(n-1)}{f(n)+1}}^k <
        \leri{1+\frac{f(n-1)}{f(n)+1}}^{n-1}=a_n,
    }
    it follows that $(b_n)_{n=1}^\infty$ is bounded as well.
    
     $(\Leftarrow)$ Suppose that $(b_n)_{n=1}^\infty$ is bounded above by some $C>0$. Since $f$ is non-decreasing, for all $n \in \N$, we have $\frac{nf(n-1)}{f(n) + 1}\le\frac{(n-1)f(n-1)}{f(n)+1}+1=b_n+1<C+1$. 
    It follows that for any $n \in \N$ with $n>1$, we have 
        \begin{gather*}
        a_n=\leri{1+\frac{f(n-1)}{f(n)+1}}^{n-1}
            <
            \leri{1 + \frac{f(n -1)}{f(n) + 1}}^n=
            \sum_{k=0}^{n} \frac{n(n-1)\cdots(n-k+1)}{k!} \leri{\frac{f(n-1)}{f(n) + 1}}^k\\
            \hspace{-5mm}<
            \sum_{k=0}^n\frac{n^k}{k!}\leri{\frac{f(n-1)}{f(n) + 1}}^k=\sum_{k=0}^n\frac{1}{k!}\leri{\frac{nf(n-1)}{f(n) + 1}}^k<\sum_{k=0}^n\frac{(C+1)^k}{k!}<\sum_{k=0}^\infty\frac{(C+1)^k}{k!}=e^{C+1}.
        \end{gather*}   
    Thus, $(a_n)_{n=1}^\infty$ is bounded above by $e^{C+1}$.
    \end{Proof}

\begin{remark}
    Let $k\in\N$, and define
        \[ F_n = \{p_1^{c_1} \cdots p_n^{c_n} : 0 \leq c_1,\ldots,c_n \leq n^k \}, \quad n=1,2,3,... \]
    Then Theorem \ref{thm 5.3} shows that $(F_n)_{n=1}^\infty$ is not a tempered Følner sequence in $(\N,\cdot)$. On the other hand, the sequence $(G_n)_{n=1}^\infty$ defined by
        \[ G_n = \{p_1^{c_1} \cdots p_n^{c_n} : 0 \leq c_i \leq n^n \text{ for } 1 \leq i \leq n \}  \]
    is a tempered Følner sequence in $(\N,\cdot)$.
\end{remark}

\subsection{Tempered Sequences in the Heisenberg group}\label{Sec 5.3}
        
\begin{Theorem} \label{heisenberg}
    Let $H$ denote the Heisenberg group over $\Z$:
        \[ H = \left\{ \begin{pmatrix}
                            1 & a & c \\
                            0 & 1 & b \\
                            0 & 0 & 1 
                        \end{pmatrix} : a,b,c \in \Z \right\}. \]
    For $n \in \N$, let
        \[ F_n = \left\{ \begin{pmatrix}
                            1 & a & c \\
                            0 & 1 & b \\
                            0 & 0 & 1 
                        \end{pmatrix} : a,b \in \{-n, \ldots, n\} \text{ and } c \in \{-n^2, \ldots, n^2\} \right\} .\]
    Then $(F_n)_{n=1}^\infty$ is a tempered Følner sequence in $H$. 
\end{Theorem}

    \begin{Proof}
        In what follows, we will write $(a,b,c)$ to denote the matrix
            \[ \begin{pmatrix}
                            1 & a & c \\
                            0 & 1 & b \\
                            0 & 0 & 1 
                        \end{pmatrix}. \]
        Then for $n \in \N$, we have $ F_n = \{ (a,b,c): a,b \in \{-n, \ldots, n\} \text{ and } c \in \{-n^2, \ldots, n^2\} \}$. It is a well-known fact that $(F_n)_{n=1}^\oo$ is a Følner sequence in $H$, so we proceed to proving that $(F_n)_{n=1}^\oo$ is tempered.
    For $(x,y,z) \in H$, we have $(x,y,z)^{-1} = (-x,-y,xy - z)$, so for $n\in\N$ with $n>1$, we see that
        \begin{align*}
        F_{n-1}^{-1} &= \{(-x,-y,xy-z): -(n-1) \leq x,y \leq n-1 \text{ and } -(n-1)^2 \leq z \leq (n-1)^2\} \\
        &= \{ (x,y,z): -(n-1) \leq x,y \leq n-1 \text{ and } -(n-1)^2 + xy \leq z \leq (n-1)^2 + xy \}.
        \end{align*}
    For $(x,y,z),(a,b,c) \in H$, one has $(x,y,z) \cdot (a,b,c) = (x + a, y + b, z + xb + c)$. Therefore, for $n > 1$,
        \begin{align*} 
            |F_{n-1}^{-1} F_n | &= \left| \left\{(x + a, y + b, z + xb + c) : 
                \begin{array}{l}
                    -(n-1) \leq x,y \leq n-1 \\
                    -(n-1)^2 + xy \leq z \leq (n-1)^2 + xy \\
                    -n \leq a,b \leq n \\
                    -n^2 \leq c \leq n^2 
                \end{array} \right\} \right| \\
            & \leq \left|\left\{x + a: \begin{array}{l}
                    -(n-1) \leq x \leq n-1 \\
                    -n \leq a \leq n \\
                \end{array}\right\}\right| \cdot \left|\left\{ y+b: \begin{array}{l}
                    -(n-1) \leq y \leq n-1 \\
                    -n \leq b \leq n \\
                \end{array} \right\}\right| \\
                & \hphantom{hi} \quad \quad \cdot \left| \left\{z + xb + c: 
                    \begin{array}{l}
                    -(n-1) \leq x,y \leq n-1 \\
                    -(n-1)^2 + xy \leq z \leq (n-1)^2 + xy \\
                    -n \leq b \leq n \\
                    -n^2 \leq c \leq n^2 
                \end{array}\right\} \right| \\
            &  \leq (2(2n-1) + 1)^2(4(n-1)^2  + 2n(n-1) + 2n^2 + 1).
        \end{align*}
    It follows that, for all $n >1$,
        \begin{align} \label{label}
            \frac{|\bigcup_{k < n} F_k^{-1}F_n|}{|F_n|} \leq \frac{(2(2n-1) + 1)^2(4(n-1)^2  + 2n(n-1) + 2n^2 + 1)}{(2n+1)^2(2n^2 + 1)}.
        \end{align} 
    Let $f(n)$ denote the expression on the right-hand side of formula (\ref{label}). Note that $\lim_{n\to\infty}f(n)$ exists and is finite. Let $C=\sup_{n\in\N}f(n)$. It follows from formula (\ref{label}) that  $\abs{\bigcup_{k < n} F_k^{-1}F_n}\le C\abs{F_n}$ for all $n\in\N$. Thus, $(F_n)_{n=1}^\infty$ is tempered. 
    \end{Proof}

\begin{remark}
    A similar argument shows that the sequence
    \equat{
            F_n=\set{(a,b,c):a,b\in\set{-p(n),...,p(n)}, c\in\set{-q(n),...,q(n)}},\quad n=1,2,3,...
    }
    forms a left tempered Følner sequence in the Heisenberg group $H$ whenever $p(x), q(x) \in \Z[x]$ are positive polynomials with $\deg(q)\ge 2\deg(p)$. 
\end{remark}
\subsection{Tempered sequences in Locally Finite Groups}\label{Sec 5.4}
        Suppose that $G$ is a countably infinite, locally finite group. Let $g_1,g_2,g_3,\ldots$ be an enumeration of its elements and for all $n \in\N,$ let $G_n$ be the (finite) subgroup generated by the elements $g_1,\ldots,g_n$. Then $(G_n)_{n=1}^\infty$ is a tempered Følner sequence in $G$.
        The next two examples give instances of this idea being applied to locally finite non-commutative groups. 

    \begin{example}
    Let $S$ denote the group of all finite permutations of $\N$: 
    \[S = \{\sigma: \N \rightarrow \N  , \,\sigma(n) = n \text{ for all but finitely many } n \}.\] 
    For $n \in \N$ the symmetric group $S_n$ is a subgroup of $S$ and, in fact, $S = \bigcup_{n=1}^\infty S_n$. If $\sigma \in S$, then $\sigma \in S_m$ for some $m$, so for $n > m$ we have $|S_n \cap \sigma S_n| = |S_n|$. It follows that $|S_n \cap \sigma S_n| / |S_n| \rightarrow 1$, hence $(S_n)_{n=1}^\infty$ is a Følner sequence. For any $n\in\N$, since $S_n$ is a group, we have $S_{n-1}^{-1}S_n \subseteq S_n$, so $(S_n)_{n=1}^\infty$ is tempered. 

    \end{example}

    \begin{example}
        Suppose that $\mathbb{F}_q$ is a finite field and define $H$ by
            \begin{equation}
                H = \left\{ \begin{pmatrix}
                            1 & f & h \\
                            0 & 1 & g \\
                            0 & 0 & 1 
                        \end{pmatrix} : f,g,h \in \mathbb{F}_q[x] \right\}.
            \end{equation}
        For $n \in \N$, define 
            \[ F_n = \left\{ \begin{pmatrix} 1 & f & h \\
            0 & 1 & g \\
            0 & 0 & 1 \end{pmatrix}: \  \deg(f), \deg(g) \leq n \text{ and } \deg(h) \leq 2n \right\}. \]
        Then $F_1 \sseq F_2 \sseq \cdots$ is an increasing sequence of subgroups of $H$ (which is locally finite) and $\bigcup_{n=1}^\oo F_n = H$. It follows that $(F_n)_{n=1}^\infty$ is a tempered Følner sequence in $H$. 
        
    \end{example}

\newpage

\bibliographystyle{alpha}
\bibliography{Arxiv}






\end{document}